\documentclass[aop,seceqn,MSNbibl,citesort,dvips]{arximspdf}

\doi{10.1214/10-AOP574}
\volume{39}
\issue{3}
\pubyear{2011}
\firstpage{939}
\lastpage{984}

\makeatletter

\newtheorem{theorem}{Theorem}[section]
\newtheorem{lemma}[theorem]{Lemma}
\newproclaim{condition}[theorem]{Condition}
\newtheorem{proposition}[theorem]{Proposition}
\newproclaim{remark}[theorem]{Remark}
\newproclaim{example}[theorem]{Example}
\mathchardef\varLambda="0103
\makeatother

\begin{document}
\begin{frontmatter}

\title{Poisson representations of branching
Markov and measure-valued branching processes}
\runtitle{Poisson representations}

\begin{aug}
\author[A]{\fnms{Thomas G.} \snm{Kurtz}\corref{}\thanksref{t1}\ead
[label=e1]{kurtz@math.wisc.edu}} and
\author[B]{\fnms{Eliane R.} \snm{Rodrigues}\thanksref{t2}\ead
[label=e2]{eliane@math.unam.mx}\ead[label=e3]{eliane@math.cinvestav.mx}}
\runauthor{T. G. Kurtz and E. R. Rodrigues}
\affiliation{University of Wisconsin, Madison and UNAM}
\address[A]{Departments of Mathematics\\
and Statistics\\
University of Wisconsin, Madison\\
480 Lincoln Drive\\
Madison, Wisconsin 53706-1388\\
USA\\
\printead{e1}} %adresu isvedimo komanda gale!
\address[B]{Instituto de Matem\'aticas\\
UNAM\\
M\'exico, DF 04510\\
Mexico \\
\printead{e2}\\
\phantom{E-mail: }\printead*{e3}}
\end{aug}

\thankstext{t1}{Supported in part by
NSF Grants DMS-05-03983 and DMS-08-05793.}

\thankstext{t2}{Supported in part by DGAPA-UNAM, Mexico, Grant 968SFA.}

% HISTORY:
\received{\smonth{1} \syear{2008}}
\revised{\smonth{11} \syear{2009}}

% ABSTRACT
%
\begin{abstract}
Representations of branching Markov processes and their
measure-valued limits in terms of
countable systems of particles are constructed for
models with spatially varying birth and death rates.
Each particle has a location and a ``level,'' but unlike
earlier constructions, the levels change with time. In
fact, death of a particle occurs only when the level of
the particle crosses a specified level $r$, or for the limiting
models, hits infinity. For branching Markov processes,
at each time $t$, conditioned on the state of the process,
the levels are independent and uniformly distributed on
$[0,r]$.
For the limiting measure-valued process, at each time $t$, the joint
distribution of locations and levels is conditionally
Poisson distributed with mean measure $K(t)\times\varLambda$,
where $\varLambda$ denotes Lebesgue measure, and $K$ is the
desired measure-valued process.

The representation simplifies or gives alternative proofs
for a variety of calculations and results
including conditioning on extinction or
nonextinction, Harris's convergence theorem for
supercritical branching processes, and diffusion
approximations for processes in random environments.
\end{abstract}

% KEYWORDS
%

\end{frontmatter}

%s1 ###
\section{Introduction}
Measure-valued processes arise naturally as infinite
system limits of
empirical measures of finite particle systems. A
number of approaches have been developed which
preserve distinct particles in the limit and which give a
representation of the measure-valued process as a
transformation of the limiting infinite particle system.
Most of these representations
\cite{DK96,DK99b,DK99,KX99} have exploited properties
of exchangeable sequences, identifying the state of the
measure-valued process as a multiple of the de Finetti
measure of the sequence, or, as in \cite{KX99}, as a
transformation of the de Finetti measure of a sequence
that gives both a location and a mass for the distinct
particles.

The primary limitation of the representations given for
measure-valued branching processes in \cite{DK99} is
that the branching rates must be independent of particle
location. This restriction was relaxed in \cite{K00} for
critical and subcritical branching, but that representation
does not seem to provide much useful insight and would
be difficult to extend to supercritical processes. A
second limitation of the approach in \cite{DK99} is that,
at least without major additional effort,
it applies only to models in which the states are finite
measures.

In the present paper, we give another construction
which, although similar to that of \cite{K00}, applies
immediately to both subcritical and supercritical
processes (as well as processes that are subcritical in
some locations and supercritical in others). The
construction also applies immediately to models with
infinite mass. In addition,
the new construction seems to be a much more effective
tool for analyzing the measure-valued processes obtained.

We introduce the basic ideas of the construction in
Section \ref{simpex}, giving results for population models without
location or type. As in the earlier work, the
justification for the representation is a consequence of a
Markov mapping theorem, Theorem~\ref{mf}. The Feller
diffusion approximation for nearly critical branching
processes is obtained as a consequence of the
construction. Section
\ref{sectpois} gives the construction for the general
branching Markov process and the Dawson--Watanabe
superprocess limit. Section \ref{secteg} gives a variety
of applications and extensions, including conditioning on
nonextinction, models with heavy-tailed offspring
distributions and processes in random environments.
The \hyperref[app]{Appendix} contains background material and a number
of technical lemmas.

%s2 ###
\section{Simple examples}\label{simpex}

In this section, we give particle representations for pure
death and continuous time Markov branching processes
and illustrate how the infinite system limit can be
derived immediately. The main point is to introduce the
notion of the \textit{level} of a particle in the simplest possible
settings.

%s2.1 ###
\subsection{Pure death processes}
For $r>0$, let $\xi_1(0),\ldots,\xi_{n_0}(0)$ be independent random variables,
uniformly distributed on $[0,r]$.
For $b>0$, let
\[
\dot{\xi}_i(t)=b\xi_i(t),
\]
so $\xi_i(t)=\xi_i(0)e^{bt}$, and define $N(t)=\#\{i\dvtx\xi_i(t)<r\}$ and
$U(t)=(U_1(t),\ldots,\break U_{N(t)})$, where the $U_j(t)$ are the values of
the $\xi_i(t)$ that are less than $r$. The $U_i(t)$ will be referred
to as the levels of the particles. The level of a particle
being below $r$ means that the particle is ``alive,'' and as
soon as its level reaches $r$ the particle ``dies.'' Note
that $N(t)$ is the number of particles ``alive'' in the
system at time $t$.

Let $f(u,n)=\prod_{i=1}^ng(u_i)$, where
$0\leq g\leq1$, $g$ is continuously differentiable and $g(u_i)=
1$
for $u_i>r$. [The ``$n$'' in $f(u,n)$ is, of course, redundant, but it will
help clarify some of the later calculations.] Then
\[
\frac d{dt}f(U(t),N(t))=Af(U(t),N(t)),
\]
where
\[
Af(u,n)=f(u,n)\sum_{i=1}^nbu_ig'(u_i)/g(u_i).
\]
Note that $Af(u,n)$ may also be written as
\[
Af(u,n)=\sum_{i=1}^nbu_ig'(u_i)\prod_{j\neq i}g(u_j).
\]
Hence, even if $g(u_i)=0$ for some $i$, the expression for
$Af(u,n)$ still makes sense.

Let $\alpha_r(n,du)$ be the joint distribution of $n$ i.i.d. uniform
$[0,r]$ random variables. Setting $e^{-\lambda_g}=r^{-1}\int_0^r
g(z)\,dz$,
$\widehat{f}(n)=\int f(u,n)\alpha_r(n,du)=e^{-\lambda_gn}$ and
\begin{eqnarray*}
\int Af(u,n)\alpha_r(n,du)&=&ne^{-\lambda_g(n-1)}br^{-1}\int_0^rz
g'(z)\,dz\\
&=&bne^{-\lambda_g(n-1)}r^{-1}\int_0^r\bigl(g(r)-g(z)\bigr)\,dz\\
&=&bne^{-\lambda_g(n-1)}(1-e^{-\lambda_g})\\
&=&C\widehat{f}(n),
\end{eqnarray*}
where
\[
C\widehat{f}(n)=bn[\widehat{f}(n-1)-\widehat{f}(n)],
\]
that is, the generator of a linear death process. Of
course, the conditional distribution of $U(t)$ given $N(t)$
is just $\alpha_r(N(t),\cdot)$.

Let $\mathcal{F}_t=\sigma(U(s)\dvtx s\leq t)$ and $\mathcal{F}^N_t=\sigma(
N(s)\dvtx s\leq t)$. Then
trivially,
\[
f(U(0),N(0))=f(U(t),N(t))-\int_0^tAf(U(s),N(s))\,ds
\]
is an $\{\mathcal{F}_t\}$-martingale, and Lemma \ref{mglem} implies
\begin{eqnarray*}
&&E[f(U(t),N(t))|\mathcal{F}^N_t]-\int_0^tE[Af(U(s),N(s))|\mathcal{F}^N_s
]\,ds\\
&&\qquad=\widehat{f}(N(t))-\int_0^tC\widehat{f}(N(s))\,ds
\end{eqnarray*}
is a $\{\mathcal{F}^N_t\}$-martingale. Consequently, $N$ is a solution of
the martingale problem for $C$ and hence is a linear death
process. Of course, this observation follows immediately
from the fact that $\tau_i^r$ defined by $U_i(0)e^{b\tau_i^r}=r$ is
exponentially distributed with parameter $b$, but the
martingale argument illustrates a procedure that works
much more generally.

%s2.2 ###
\subsection{A simple branching process}\label{simpbp}
For $f$ and $g$ as above, $a>0$, $r>0$ and $-\infty<b\leq ra$, define the
generator
%
%e2.1 ###
\begin{eqnarray}\label{lbgen}
A_rf(u,n)&=&f(u,n)\sum_{i=1}^n2a\int_{u_i}^r\bigl(g(v)-1
\bigr)\,dv\nonumber\\[-8pt]\\[-8pt]
&&{}+f(u,n)\sum_{i=1}^n(au_i^2-bu_i)\frac{g'(u_i)}{g(u_i)}.\nonumber
\end{eqnarray}
We refer to $u_i$ as the \textit{level} of the $i$th particle, and as
in the pure death example, a~particle ``dies'' when its
level reaches $r$. The process with generator (\ref{lbgen})
has the following properties. The particle levels satisfy
\[
\dot{U}_i(t)=aU_i^2(t)-bU_i(t),
\]
and a particle with level $z$ gives birth at rate $2a(r-z)$ to
a particle whose initial level is uniformly distributed
between $z$ and $r$. Uniqueness for the martingale
problems for $A_r$ follows by first checking uniqueness
for the operator $D$ given by the second term alone, and then
showing that uniqueness holds up to the first time
that the process includes $n$ particles by observing that
$A_r$ truncated at $n$ is a bounded perturbation of $D$.
Finally, the first hitting time of $n$ goes to infinity as
$n\rightarrow\infty$ (cf. Problem 28 in Section 4.11 of \cite{EK86}).

As before,
$\widehat{f}(n)=\int f(u,n)\alpha_r(n,du)=e^{-\lambda_gn}$. To calculate
$\int A_r(fu,n)\times\alpha_r(n,du)$, observe that
\[
r^{-1}2a\int_0^rg(z)\int_z^r\bigl(g(v)-1\bigr)\,dv\,dz=are^{-2\lambda_g}-2ar^{
-1}\int_0^rg(z)(r-z)\,dz
\]
and
\begin{eqnarray*}
r^{-1}\int_0^r(az^2-bz)g'(z)\,dz&=&-r^{-1}\int_0^r(2az-b)\bigl(g(z)-1\bigr)\,dz\\
&=&-2ar^{-1}\int_0^rzg(z)\,dz+ar+b(e^{-\lambda_g}-1).
\end{eqnarray*}
Then
\begin{eqnarray*}
\int Af(u,n)\alpha_r(n,du)&=&ne^{-\lambda_g(n-1)}\bigl(are^{-2\lambda_
g}-2are^{-\lambda_g}+ar+b(e^{-\lambda_g}-1)\bigr)\\
&=&C\widehat{f}(n),
\end{eqnarray*}
where
%
%e2.2 ###
\begin{equation}\label{brproj}
C\widehat{f}(n)=ran\bigl(\widehat{f}(n+1)-\widehat{f}(n)\bigr)+(ra-b
)n\bigl(\widehat{f}(n-1)-\widehat{f}(n)\bigr)
\end{equation}
is the generator of a branching process.

Unlike the linear death example, it is not immediately
obvious that\break $\alpha_r(N(t),\cdot)$ is the conditional distribution of
$U(t)$ given $\mathcal{F}_t^N=\sigma(N(s)\dvtx s\leq t)$; however, Theorem
\ref{mf} and
the fact that a solution of the martingale problem for $C$
starting from $N(0)$ exists gives the existence of a
solution $(U(t),N(t))$ of the martingale problem for $A_{r}$ such
that for all $t\geq0$, $\alpha_r(N(t),\cdot)$ is the conditional
distribution of $U(t)$
given $\mathcal{F}^N_t$. To apply Theorem \ref{mf}, take $\psi(u,n
)=n$ in
(\ref{opest}) and assume $E[N(0)]<\infty$. Any solution of the
martingale problem for $C$ with $E[N(0)]<\infty$ will satisfy
$E[N(t)]=E[N(0)]e^{bt}$, for all $t\geq0$. The moment assumption can be
eliminated by conditioning.

We conclude that for any distribution for $N(0)$, there
is a solution $(U,N)$ of the martingale
problem for $A$ such that $N$ is a solution of the
martingale problem for~$C$, that is, $N$ is a linear birth
and death process with birth rate $ar$ and death rate
$ar-b$. Uniqueness holds for the martingale problem for
$A$, so for any solution of the martingale problem for $A$
satisfying $P\{U(0)\in\Gamma|N(0)\}=\alpha_r(N(0),\Gamma)$, we have
that $
N$
is a solution of the martingale problem for $C$.

This representation can be used to do simple
calculations. For example, let $U_{*}(0)$ be the minimum of
$U_1(0),\ldots,U_{N(0)}$. Then for all $t$, all levels are
above
\[
U_{*}(t)=\frac{U_{*}(0)e^{-bt}}{1-(a/b)U_{*}(0)(1-e^{-bt})}
.
\]
Let $\tau=\inf\{t\dvtx N(t)=0\}$. Then if $\tau$ is finite, $U_{*}(\tau
)=r$. In
particular, if $N(0)=n$, then
\begin{eqnarray*}
P\{\tau>t\}&=&P\{U_{*}(t)<r\}=P\biggl\{U_{*}(0)<\frac r{e^{-bt}-({r
a}/b)(e^{-bt}-1)}\biggr\}\\
&=&1-\biggl(e^{-bt}-\frac{ra}b(e^{-bt}-1)\biggr)^{-n}.
\end{eqnarray*}
Note that the assumption that $b\leq ra$ ensures that
$e^{-bt}-\frac{ra}b(e^{-bt}-1)\geq1$.

In the branching process, the average lifetime of an
individual is $(ar-b)^{-1}$ which will be small if $r$ is large.
Consequently, it is important to note that the levels do
not represent single individuals in the branching process
but whole lines of descent. For example, at least in the
critical or subcritical case, the individual with level
$U_{*}(0)$ at time zero is the individual whose line of
descent lasts longer than that of any other individual
alive at time zero.

%s2.2.1 ###
\subsubsection{Conditioning on nonextinction}\label{cndne}
If $b\leq 0$, then
conditioning on non\-extinction, that is, conditioning on
$\tau>t$ and letting $t\rightarrow\infty$, is equivalent to
conditioning on
$U_{*}(0)=0$. Conditioned on $U_{*}(0)=0$, $U_1(0),\ldots,U_{N(
0)}(0)$
include $N(0)-1$ independent, uniform $[0,r]$ random
variables and one that equals zero. If one of the initial
levels is zero, then the solution of the martingale
problems for $A_r$ gives a solution for
\begin{eqnarray*}
A_r^cf(u,n)&=&f(u,n)\sum_{i=1}^{n-1}2a\int_{u_i}^r\bigl(g(v)-1\bigr)\,dv+f(u,n)
\sum_{i=1}^{n-1}(au_i^2-bu_i)\frac{g'(u_i)}{g(u_i)}\\
&&{} + f(u,n)2a\int_0^r\bigl(g(v)-1\bigr)\,dv,
\end{eqnarray*}
and taking $\alpha^c_r(n,du)$ to be the distribution of $n$
independent random variables, one of which is zero and
the others uniform $[0,r]$, we see that
$N$ is a solution of the martingale problem for
\[
C^c\widehat{f}(n)=ra(n+1)\bigl(\widehat{f}(n+1)-\widehat{f}(n)\bigr)+(ra-b)(n-1)\bigl(\widehat{
f}(n-1)-\widehat{f}(n)\bigr).
\]

%s2.2.2 ###
\subsubsection{Conditioning on extinction}\label{cnde}
If $0<b<ra$,
then conditioning on extinction is equivalent to
conditioning on $U_{*}(0)>\frac ba$. Conditioned on $U_{*}(0)>\frac
ba$,
$U_1(0),\ldots,U_{N(0)}(0)$ are independent uniform $[\frac ba,r
]$.
Defining $V_i(t)=U_i(t)-\frac ba$, $V$ is a solution of the
martingale problem for
\[
A_rf(v,n)=f(v,n)\sum_{i=1}^n2a\int_{v_i}^{r-b/a}\bigl(g(z)-1\bigr)\,dz+f(v,
n)\sum_{i=1}^n(av_i^2+bv_i)\frac{g'(v_i)}{g(v_i)},
\]
so $N$ is a solution of the martingale problem for
\[
C\widehat{f}(n)=(ra-b)n\bigl(\widehat{f}(n+1)-\widehat{f}(n)\bigr)+ran\bigl(\widehat{f}(n-
1)-\widehat{f}(n)\bigr),
\]
which is the generator of a subcritical branching
process.

%s2.2.3 ###
\subsubsection{Convergence as $t\rightarrow\infty$}\label{harthm}
Again, in the supercritical case, $0<b<ra$, if
$0<U_{*}(0)<\frac ba$, then $N(t)\rightarrow\infty$. Observe that
\[
V_{*}(\infty)=\lim_{t\rightarrow\infty}e^{bt}U_{*}(t)=\frac{U_{
*}(0)}{1-(a/b)U_{*}(0)}
\]
exists, and a similar limit will hold for any level whose
initial value is below $\frac ba$. Setting $\xi(t)=\sum\delta_{
e^{bt}U_i(t)}$, the
counting measure $\xi(t)$ converges almost surely in the
sense that
\[
\lim_{t\rightarrow\infty}\int_0^{\infty}f(u)\xi(t,du)=\lim
_{t\rightarrow
\infty}\sum_if(e^{bt}U_i(t))=\int_0^{\infty}f(u)\xi(\infty,du)
\qquad\mbox{a.s.}
\]
for each bounded, continuous, nonnegative
$f$ with compact support in
$[0,\infty)$. Let $\{\mathcal{F}^N_t\}$ be the filtration generated by $
N$. Then
as in (\ref{uniid}),
%
%e2.3 ###
\begin{equation}\label{popnrm}
E\bigl[e^{-\int_0^{re^{bt}}f(u)\xi(t,du)}|\mathcal{F}_t^
N\bigr]=e^{-F_f^tr^{-1}e^{-bt}N(t)},
\end{equation}
where
\[
F^t_f=-re^{bt}\log\biggl(1-r^{-1}e^{-bt}\int_0^{re^{bt}}\bigl(1-e^{-f(u)}\bigr)\,
du\biggr)\rightarrow\int_0^{\infty}\bigl(1-e^{-f(u)}\bigr)\,du.
\]
The left-hand side of (\ref{popnrm}) converges almost surely
by Lemma \ref{cnddc}. Consequently,
\[
W\equiv\lim_{t\rightarrow\infty}e^{-bt}N(t)
\]
exists almost surely. Note that $W>0$ if and only if
$\lim_{t\rightarrow\infty}N(t)=\infty$.

Conditioned on $W$, $\xi(\infty)$ is a
Poisson point process with intensity $r^{-1}W$, and
$V_{*}(\infty)$ is exponentially distributed with
parameter $r^{-1}W$, with the understanding that $V_{*}(\infty)=
\infty$ if
$W=0$. It follows that for $\lambda>0$,
\begin{eqnarray*}
&&
P\{r^{-1}V_{*}(\infty)>\lambda|V_{*}(\infty)<\infty\}\\
&&\qquad=E[e^{
-\lambda W}|W>0]\\
&&\qquad=P\biggl\{r^{-1}U_{*}(0)\biggl(1-\frac abU_{*}(0)\biggr)^{-1}>\lambda\Big|U_{*}(0)<\frac
ba\biggr\}.
\end{eqnarray*}
If $N(0)=1$, then $U_{*}(0)$ is uniformly distributed on $[0,r]$,
and hence $P\{W>0\}=\frac b{ra}$ and
\[
E[e^{-\lambda W}|W>0]=\frac1{1+({ra}/b)\lambda},
\]
that is, $W$ is exponentially distributed with parameter
$\frac b{ra}$. Of course, we have simply rederived a classical
result of Harris \cite{Har51}.

%s2.3 ###
\subsection{Feller diffusion approximation}\label{fdiff}
As $r\rightarrow\infty$,
$A_rf$ in (\ref{lbgen})
converges for every continuously differentiable
$g$ such that $0\leq g\leq1$ and $g(z)=1$ for
$z\geq r_g$, that is, for $f(u)=\prod_ig(u_i)$, in the limit
%
%e2.4 ###
\begin{equation}\label{felldiff}\qquad
Af(u)=f(u)\sum_i2a\int_{u_i}^{r_g}\bigl(g(v)-1\bigr)\,dv+f(u)
\sum_i(au_i^2-bu_i)\frac{g'(u_i)}{g(u_i)}.
\end{equation}
If $n/r\rightarrow y$ as $r\rightarrow\infty$, then $\alpha_r(n,\cdot
)$ converges to $\alpha(y,\cdot)$,
where $\alpha(y,\cdot)$ is the distribution of a Poisson process
$\xi_y$ on
$[0,\infty)$ with intensity $y$, in the sense that
\[
\int_{[0,r]^n}f(u)\alpha_r(n,du)\rightarrow E\bigl[e^{\int_0^{\infty}\log
g(z)\,d\xi_y(z)}\bigr]=e^{-y\int_0^{\infty}(1-g(z))\,dz}.
\]
Note that
\[
\widehat{f}(y)=\alpha f(y)=\int f(u)\alpha(y,du)=e^{-y\int_0^{\infty}
(1-g(z))\,dz}=e^{-y\beta_g},
\]
and using Lemma \ref{lemma2}
\begin{eqnarray*}
\alpha Af(y)&=&e^{-y\beta_g}\biggl(2ay\int_0^{\infty}g(z)\int_z^{\infty}
\bigl(g(v)-1\bigr)\,dv\,dz\\
&&\hspace*{62.3pt}{}+y\int_0^{\infty}(az^2-bz)g'(z)\,dz\biggr)\\
&=&e^{-y\beta_g}\biggl(2ay\int_0^{\infty}g(z)\int_z^{\infty}\bigl(g(v)-
1\bigr)\,dv\,dz\\
&&\hspace*{43.3pt}{}-y\int_0^{\infty}(2az-b)\bigl(g(z)-1\bigr)\,dz\biggr)\\
&=&e^{-y\beta_g}\biggl(2ay\int_0^{\infty}g(z)\int_z^{\infty}\bigl(g(v)-1\bigr)\,dv\,dz\\
&&\hspace*{30pt}{}-2ay\int_0^{\infty}\int_z^{\infty}\bigl(g(v)-1\bigr)\,dv\,dz\\
&&\hspace*{81.3pt}{} +by\int_0^{\infty}\bigl(g(z)-1\bigr)
\,dz\biggr)\\
&=&e^{-y\beta_g}(ay\beta_g^2-by\beta_g)\\
&=&C\widehat{f}(y),
\end{eqnarray*}
where
\[
C\widehat{f}(y)=ay\widehat{f}^{\hspace*{1pt}\prime\prime}(y)+by\widehat{f}^{\hspace*{1pt}\prime}(y).
\]

Again, we can apply Theorem \ref{mf} taking
\[
\gamma(u)=\limsup_{z\rightarrow\infty}\frac1z\sum_i\mathbf{1}_{[
0,z]}(u_i)
\]
and $\psi(u)=\sum_ie^{-u_i}$, so
\[
|Af(u)|\leq[2ar_g+\Vert g'\Vert(ar_g^2+|b|r_g)] e^{r_g}\psi(u)
\]
and\vspace*{1pt} $\int\psi(u)\alpha(y,du)=y$. If $\widetilde{Y}$ is a solution of the
martingale problem for $C$ with $E[\widetilde{Y}(0)]<\infty$, then
$E[\widetilde{Y}(t)]=e^{bt}E[\widetilde{Y}(0)]$ and the conditions of Theorem
\ref{mf}
are satisfied. Consequently, there is a solution $U$ of the
martingale problem for $A$ such that
%
%e2.5 ###
\begin{equation}\label{ycomp}
Y(t)=\limsup_{z\rightarrow\infty}\frac1z\sum_i{\mathbf1}_{
[0,z]}(U_i(t))
\end{equation}
is a solution of the martingale problem for $C$ with the
same distribution as $\widetilde{Y}$. [Note
that, with probability one, the $\limsup$ in (\ref{ycomp}) is
actually a limit.]

If $U$ is a solution of the martingale problem for $A$, then
$U^r(t)=\{U_i(t)\dvtx\break U_i(t)<r\}$ defines a solution of the
martingale problem for $A_r$. Uniqueness for $A_r$ follows
by the argument outlined in Section \ref{unmgp}, and
uniqueness for $A_r$
implies uniqueness for $A$. Since uniqueness holds for the
martingale problem for $A$, by Theorem~\ref{mf}(c),
uniqueness holds for $C$ also. In general, if $U$ is a
solution of the martingale problem for $A$ and $\sum_i\delta_{U_i
(0)}$ is
a Poisson random measure with mean measure $y\varLambda$, where
$\varLambda$ denotes Lebesgue measure, then (\ref{ycomp}) is a
solution of the martingale problem for $C$.

%s2.4 ###
\subsection{The genealogy and the number of
ancestors}\label{geneal}
For each $T>0$, there is a solution of
%
%e2.6 ###
\begin{equation}\label{popbnd}
\dot{u}_T(t)=au_T(t)^2-bu_T(t)
\end{equation}
satisfying $u_T(t)<\infty$ for $t<T$ and $\lim_{t\rightarrow T-}u_
T(t)=\infty$.
Every particle alive at time $T$ is a descendent of some
particle $U_i(t)$ alive at time $t<T$ satisfying
$U_i(t)<u_T(t)$. Note that the converse is also true. If
$U_i(t)<u_T(t)$, then $U_i(t)$ has descendants alive at
time $T$. In fact, a positive fraction of the particles
alive at time $T$ will be descendants of $U_i(t)$.

If $Y(T)>0$, then there are
infinitely many particles alive at time $T$, but since
\[
\xi(t,[0,u_T(t)))<\infty,
\]
they are
all descendants of finitely many ancestors alive at time
$t$. Note that $t\rightarrow\xi(t,[0,u_T(t)))$ is nondecreasing and
increases by jumps of $+1$. It is not possible to recover
the full genealogy just from the levels since a new
individual appearing at time $t$ with level $v$ could be the
offspring of any existing individual with level $U_i(t)<v$.
In Section \ref{sectpois}, particles will be assigned a
location (or type), and if these locations evolve in such
a way that two particles have the same location only if
one is the offspring of the other and then only at the
instant of birth, it will be possible to reconstruct the
full genealogy from the levels and locations.

%s2.5 ###
\subsection{Branching processes in random environments}
Assume that $a$ and $b$ are functions of another stochastic
process $\xi$, say an irreducible, finite Markov chain with
generator $Q$. Then, for functions of the form
$f(l,u,n)=f_0(l)f_1(u)=f_0(l)\prod_{i=1}^ng(u_i)$,
consider a scaled generator
\begin{eqnarray*}
A_rf(l,u,n)&=&rf_1(u)Qf_0(l)+f(l,u,n)\sum_{i=1}^n2a(l)\int_{u_i}^r\bigl(
g(v)-1\bigr)\,dv\\
&&{}+f(l,u,n)\sum_{i=1}^n\bigl(a(l)u_i^2-\sqrt{r}b(l)u_i\bigr)\frac{g'(u_i)}{
g(u_i)},
\end{eqnarray*}
which, as in (\ref{brproj}), corresponds to a process with
generator
\begin{eqnarray*}
C_r\widehat{f}(l,n)&=&rQ\widehat{f}(l,n)+a(l)rn\bigl(\widehat{f}(l,n+1)-\widehat{f}
(n)\bigr)\\
&&{} +\bigl(ra(l)-\sqrt{r}b(l)\bigr)n\bigl(\widehat{f}(l,n-1)-\widehat{f}(l,n)\bigr),
\end{eqnarray*}
where $\widehat{f}(l,n)=f_0(l)e^{-\lambda_gn}$. The process corresponding to
$C_r$ is a branching process in a random environment
determined by $\xi$. Writing the process corresponding to
$A_r$ as
\[
\bigl(\xi(rt),U_1(t),\ldots,U_{N_r(t)}\bigr)
\]
the process
corresponding to $C_r$ is $(\xi(rt),N_r(t))$.

Note that in this example, the levels satisfy
\[
\dot{U}_i(t)=a(\xi(rt))U_i^2(t)-\sqrt{r}b(\xi(rt))U_i(t).
\]
Let $\pi$ be the stationary distribution for $Q$, and assume
that $\sum_l\pi(l)b(l)=0$. Then, by Theorem 2.1 or
\cite{Bha82}, for example,
\[
Z^{(r)}(t)=\sqrt{r}\int_0^tb(\xi(rs))\,ds
\]
converges to a Brownian motion $Z$ with variance
parameter
\[
\sum_k\sum_l\pi(k)q_{kl}\bigl(h_0(l)-h_0(k)\bigr)^2=-2\sum_l\pi(l)h_0(l
)b(l)\equiv2\overline{c},
\]
where $h_0(l)$ is a solution of
$Qh_0(l)=b(l)$. In the limit, by Theorem 5.10 of
\cite{KP91a} (applying a truncation argument to extend
the boundedness assumption), the levels will satisfy
%
%e2.7 ###
\begin{equation}\label{leveq4}
dU_i(t)=\bigl(\overline{a}U_i(t)^2+\overline{c}U_i(t)\bigr)\,dt+\sqrt{
2\overline{c}}U_i(t)\,dW(t),
\end{equation}
where $\overline{a}=\sum\pi(l)a(l)$.

Applying ideas from \cite{Kur73}, we can obtain
convergence for the full system by considering the
asymptotic behavior of the generator. Setting
\[
h_1(l,u,n)=h_0(l)f_1(u,n)\sum_{i=1}^nu_i\frac{g'(u_i)}{g(u_i)}
,
\]
we have
\begin{eqnarray*}
&&A_r\biggl(f_1+\frac1{\sqrt{r}}h_1\biggr)(l,u,n)\\
&&\qquad =f_1(u,n)\sum_{i=1}^n2a(l)\int_{u_i}^r\bigl(g(v)-1\bigr)\,dv+f_1(u,n)
\sum_{i=1}^na(l)u_i^2\frac{g'(u_i)}{g(u_i)}\\
&&\qquad\quad{} +\frac1{\sqrt{r}}h_0(l)f_1(u,n)\Biggl(\sum_{i=1}^nu_i\frac{
g'(u_i)}{g(u_i)}\Biggr)\sum_{i=1}^n2a(l)\int_{u_i}^r\bigl(g(v)-1\bigr)\,dv\\
&&\qquad\quad{} +\frac1{\sqrt{r}}\sum_{i=1}^nh_0(l)f_1(u,n)\int_{u_i}^
rvg'(v)\,dv\\
&&\qquad\quad{} +\frac1{\sqrt{r}}h_0(l)f_1(u,n)\sum_{j=1}^n\bigl(a(l)u_j^2-\sqrt{
r}b(l)u_j\bigr)\\
&&\qquad\quad\hspace*{103.8pt}{}\times\biggl(\sum_{i\neq j}u_i\frac{g'(u_i)g'(u_j)}{g(u_i)g(u_
j)}+\frac{g'(u_j)+u_jg^{\prime\prime}(u_j)}{g(u_j)}\biggr),
\end{eqnarray*}
and passing to the limit as $r\rightarrow\infty$, $A_r(f_1+\frac
1{\sqrt{r}}h_1)$ converges
to
\begin{eqnarray*}
&&\widetilde{A}f_1(u,l)\\
&&\qquad=f_1(u)\sum_i2a(l)\int_{u_i}^{\infty}\bigl(g(v)-1\bigr)\,dv+f_1(u)\sum_
ia(l)u_i^2\frac{g'(u_i)}{g(u_i)}\\
&&\qquad\quad{} -h_0(l)b(l)f_1(u)\sum_j\biggl(\sum_{i\neq j}u_ju_i\frac{
g'(u_i)g'(u_j)}{g(u_i)g(u_j)}+\frac{u_jg'(u_j)+u_j^2g^{\prime\prime}
(u_j)}{g(u_j)}\biggr).
\end{eqnarray*}
Finally, we can find an additional perturbation $h_2$ so that
$A_r(f_1+\frac1{\sqrt r}h_1+\frac1rh_2)$ converges to
\begin{eqnarray*}
Af_1(u)
&=&f_1(u)\sum_i2\overline{a}\int_{u_i}^{\infty}\bigl(g(v)-1\bigr)\,dv+f_1(u
)\sum_i\overline{a}u_i^2\frac{g'(u_i)}{g(u_i)}\\
&&{} +\overline{c}f_1(u)\sum_j\biggl(\sum_{i\neq j}u_ju_i\frac{g'
(u_i)g'(u_j)}{g(u_i)g(u_j)}+\frac{u_jg'(u_j)+u_j^2g^{\prime\prime}
(u_j)}{g(u_j)}\biggr).
\end{eqnarray*}
This
convergence assures convergence of the finite models to
an infinite particle model. The particle birth process is
the same as in Section \ref{fdiff}, but the levels satisfy
(\ref{leveq4}) where the Brownian motion $W$ is the
same for all levels.

Let $\alpha$ and $\beta_g$ be as in Section \ref{fdiff}, and note that
\begin{eqnarray*}
\beta_g&=&\int_0^{\infty}\bigl(1-g(z)\bigr)\,dz=\int_0^{\infty}zg'(z)\,dz\\
&=&-\frac
12\int_0^{\infty}z^2g^{\prime\prime}(z)\,dz.
\end{eqnarray*}
We have from Lemma \ref{lemma2}
\begin{eqnarray*}
\alpha Af(y)&=&e^{-y\beta_g}\biggl(2\overline{a}y\int_0^{\infty}g(z)\int_
z^{\infty}\bigl(g(v)-1\bigr)\,dv\,dz\\
&&\hspace*{30.1pt}{} + y\int_0^{\infty}(\overline{a}z^2+\overline{c}z)g'(
z)\,dz\\
&&\hspace*{30.1pt}{} + \overline{c}y^2\biggl(\int_0^{\infty}zg'(z)\,dz\biggr)^2+\overline{
c}y\int_0^{\infty}z^2g^{\prime\prime}(z)\,dz\biggr)\\
&=&e^{-y\beta_g}\bigl((\overline{a}y+\overline{c}y^2)\beta_g^2-\overline{c}y\beta_
g\bigr)\\
&=&C\widehat{f}(y),
\end{eqnarray*}
where
\[
C\widehat{f}(y)=(\overline{a}y+\overline{c}y^2)\widehat{f}^{\hspace*{1pt}\prime\prime}(y)+
\overline{c}y\widehat{f}^{\hspace*{1pt}\prime}(y),
\]
which identifies the diffusion limit for $r^{-1}N_r$.

Theorem \ref{mf}  can be extended to cover models
with non-Markovian environments, that is, the process $\xi$
is specified directly rather than through a generator.
The diffusion limit is then obtained by verifying
convergence of the level processes and applying Theorem
\ref{pmlim}.

For early work on diffusion approximations
for branching processes in
random environments, see \cite{Hell81,Keid75,Kur78} and
also \cite{EK86}, Section 9.3.

%s3 ###
\section{Representations of measure-valued branching
processes}\label{sectpois}

%s3.1 ###
\subsection{Branching Markov processes}\label{sec-bran} We now
consider particles with both a level $u_i$ and a location $x_i$
in a complete, separable metric space $E$. Since the
indexing of the particles is not important, we identify a
state $(x,u,n)$ of our process with the counting measure
$\mu_{(x,u)}=\sum_{i=1}^n\delta_{(x_i,u_i)}$. Let
\[
f(x,u,n)=\prod_{i=1}^ng(x_i,u_i)=e^{\int\log g\,d\mu_{(x,u)}},
\]
where $g\dvtx E\times[0,\infty)\rightarrow(0,1]$. We assume that as a
function of $x$, $g$ is in the domain $\mathcal{D}(B)$ of the generator
of a Markov process in $E$, $g$ is continuously differentiable
in $u$ and $g(x,u)=1$ for $u\geq r$.
We set
%
%e3.1 ###
\begin{eqnarray}\label{fingen}
A_rf(x,u,n)&=&f(x,u,n)\sum_{i=1}^n\frac{Bg(x_i,u_i)}{g(x_i,u_i)}\nonumber\\
&&{}+f
(x,u,n)\sum_{i=1}^n2a(x_i)\int_{u_i}^r\bigl(g(x_i,v)-1\bigr)\,dv\\
&&{} + f(x,u,n)\sum_{i=1}^n\bigl(a(x_i)u_i^2-b(x_i)u_i\bigr)\,\frac{\partial_{u_i}
g(x_i,u_i)}{g(x_i,u_i)}.\nonumber
\end{eqnarray}
Each particle has a location $X_i(t)$ in $E$ and a level $U_i(t)$ in
$[0,r]$. The locations evolve independently as Markov
processes with generator $B$; the levels satisfy
%
%e3.2 ###
\begin{equation}\label{leveq}
\dot{U}_i(t)=a(X_i(t))U_i^2(t)-b(X_i(t))U_i(t);
\end{equation}
particles give birth at rates $2a(X_i(t))(r-U_i(t))$; the initial
location of a new particle is the location of the parent at
the time of birth;
and
the initial level is uniformly distributed on $[U_i(t),r]$.
Particles that reach level $r$ die. Setting
$e^{-\lambda_g(x_i)}=\widehat{g}(x_i)=r^{-1}\int_0^rg(x_i,z)\,dz$ and $
\widehat{f}(x,n)=\prod_{i=1}^n\widehat{g}(x_i)=e^{-\sum_{i=1}^n\lambda_
g(x_i)}$,
calculating as in Section \ref{simpbp}, we have
\begin{eqnarray*}
C_r\widehat{f}(x,n)&=&\sum_{i=1}^nB_{x_i}\widehat{f}(x,n)+\sum_{i=1}^nra(
x_i)\bigl(\widehat{f}\bigl(b(x|x_i),n+1\bigr)-\widehat{f}(x,n)\bigr)\\
&&{} + \sum_{i=1}^n\bigl(ra(x_i)-b(x_i)\bigr)\bigl(\widehat{f}\bigl(d(x|x_i),n-1\bigr)-\widehat{f}(x,n
)\bigr),
\end{eqnarray*}
where $B_{x_i}$ is the generator $B$ applied to $\widehat{f}(x,n)$ as a
function of $x_i$, $b(x|x_i)$ is the collection of $n+1$ particles in
$E$ obtained from $x$ by adding a copy of the $i$th particle $x_i$,
and $d(x|x_i)$ is the collection of $n-1$ particles obtained
from $x$ by deleting the $i$th particle, that is, if $\mu_x$
denotes $\sum_{i=1}^n\delta_{x_i}$, then for $z\in E$,
\[
\mu_{b(x|z)}=\delta_z+\sum_{i=1}^n\delta_{x_i},\qquad
\mu_{d(x|x_j)}=\sum_{i=1}^n\delta_{x_i}-\delta_{x_j}.
\]

If $ra(z)-b(z)\geq0$ for all $z\in E$, then
$C$ is the generator of a branching Markov process with
particle motion determined by $B$, the birth rate for a
particle at $z\in E$ given by $ra(z)$ and the death rate
given by $ra(z)-b(z)$.

With Theorem \ref{mf} in mind, we make the following
assumptions on $B$, $a$, $b$ and~$r$. $C(E)$ is the space of
continuous functions on $E$, $\overline{C}(E)$ the space of bounded
continuous functions on $E$ and $M(E)$ the space of Borel
measurable functions on~$E$.
\begin{condition}\label{bmpc}
\begin{enumerate}[(iii)]
\item[(i)]$B\subset\overline{C}(E)\times C(E)$, $\mathcal{D}(B)$ is closed
under multiplication
and is separating.

\item[(ii)]$f\in\mathcal{D}(A)$ satisfies $f(x,u,n)=\prod_{i=1}^ng(x_
i,u_i)$, where
$g(z,v)=\prod_{l=1}^m(1-g_1^l(z)g_2^l(v))$ for $g_1^l\in\mathcal{D}
(B)$, $g_2^l$ differentiable\vspace*{1pt}
with support in $[0,r]$ and $0\leq g_1^l(z)g_2^l(v)\leq\rho_g<1$ for
all $
l$,
$z$ and $v$.

\item[(iii)] There exists $\psi_B\in C(E)$, $\psi_B\geq0$ and constants
$c_g\geq0$, for each $g$ in (ii), such that
\[
{\sup_u}|Bg(x,u)|\leq c_g\psi_B(x), \qquad x\in E.
\]

\item[(iv)]Defining $B_0=\{(g,(\psi_B\vee1)^{-1}Bg)\dvtx g\in\mathcal{D}
(B)\}$, $B_0$ is graph
separable (see Section \ref{mpsect}).

\item[(v)]$a,b\in M(E)$, $a\geq0$, $r>0$ and $ra-b\geq0$.
\end{enumerate}
\end{condition}

We have the following generalization of the results of
Section \ref{simpbp}.
\begin{theorem}\label{finbmp}
Assume Condition \ref{bmpc}. For $x\in E^n$ and $u\in[0,\infty
)^n$,
let
\[
\psi(x,u)=1+\sum_{i=1}^n\bigl(\psi_B(x_i)+a(x_i)+|b(x_i)|\bigr)e^{-u_i}
\]
and
\[
\widetilde{\psi}(x)=1+\sum_{i=1}^n\bigl(\psi_B(x_i)+a(x_i)+|b(x_i)|\bigr)(1-
e^{-r}).
\]

If $X$ is a solution of the martingale
problem for $C$ satisfying
%
%e3.3 ###
\begin{equation}\label{momfb}
E\biggl[\int_0^t\widetilde{\psi}(X(s))\,ds\biggr]<\infty
\qquad\mbox{for all }
t\geq0,
\end{equation}
then there is a solution $(\widetilde{X},\widetilde{U})$ of the
martingale problem for $A_r$ such that $X$ and $\widetilde{X}$ have the
same distribution.
\end{theorem}
\begin{remark}
For many models, $\psi_B$, $a$ and $b$ will be uniformly
bounded, and the moment conditions (\ref{momfb}) will
hold as long as $E[X(0)]<\infty$.
\end{remark}
\begin{pf*}{Proof of Theorem \ref{finbmp}}
Note that
\[
|A_rf(x,u,n)|\leq\bigl(2r+(1+r^2+r)d_g\bigr)e^r\psi(x,u),
\]
where $d_g$ depends on the $g_1^l$, $g_2^l$, $\partial_vg_2^l$ and $
B\prod g_1^{l_k}$ for
all choices of $\{l_1,\ldots, l_j\}\subset\{1,\ldots,m\}$.
The result then follows by application of Theorem \ref{mf}.
\end{pf*}

Theorem \ref{finbmp} applies to finite branching Markov
processes. Similar results also hold for locally finite
processes.
\begin{theorem}\label{finbmp2}
In addition to Condition \ref{bmpc}, assume that
\[
\int_0^{\infty}|g(x,u)-1|\,du+\sup_u(u+u^2)\,\partial_u g(x,u)\leq c_
g\psi_B(x).
\]
For $x\in E^{\infty}$ and $u\in[0,\infty)^{\infty}$,
let
\[
\psi(x,u)=1+\sum_{i=1}^{\infty}\psi_B(x_i)\bigl(1+a(x_i)+|b(x_i)|\bigr)e^{
-u_i}
\]
and
\[
\widetilde{\psi}(x)=1+\sum_{i=1}^{\infty}\psi_B(x_i)\bigl(1+a(x_i)+|b(x_
i)|\bigr)(1-e^{-r}).
\]

If $X$ is a solution of the martingale
problem for $C$ satisfying
%
%e3.4 ###
\begin{equation}\label{momfb2}
E\biggl[\int_0^t\widetilde{\psi}(X(s))\,ds\biggr]<\infty\qquad
\mbox{for all }
t\geq0,
\end{equation}
then there is a solution $(\widetilde{X},\widetilde{U})$ of the
martingale problem for $A_r$ such that $X$ and $\widetilde{X}$ have the
same distribution.
\end{theorem}
\begin{pf}
Note that
\[
A_rf(x,u,n)\leq2c_g e^r\psi(x,u).
\]
The result then follows by application of Theorem \ref{mf}.
\end{pf}
\begin{example}
Suppose that $ a$ and $b$ are bounded, and $B$ is a diffusion
operator with bounded drift and diffusion coefficients.
Then we can take $\mathcal{D}(B)$ to be the collection of
nonnegative $C^2$-functions with compact support and
$\psi_B(z)=\frac1{(1+|z|^2)^{\beta}}$, for $\beta>0$. If
$E[\sum_i\psi_B(X_i(0))]<\infty$, then there exists a solution of the
martingale problems for $C$ satisfying
$\sup_{s\leq t}E[\sum_i\psi_B(X_i(s))]<\infty$ and hence $E[\int_
0^t\widetilde{\psi}(X(s))\,ds]<\infty$.
\end{example}

%s3.2 ###
\subsection{Basic limit theorem}\label{rlim}
As in Section \ref{fdiff}, if $r\rightarrow\infty$,
$Af$ given by
(\ref{fingen}) becomes
%
%e3.5 ###
\begin{eqnarray}\label{eq33}
Af(x,u)&=&f(x,u)\sum_i\frac{Bg(x_i,u_i)}{g(x_i,u_i)}\nonumber\\
&&{}+f(x,u)\sum_i2
a(x_i)\int_{u_i}^{r_g}\bigl(g(x_i,v)-1\bigr)\,dv\\
&&{}+f(x,u)\sum_i\bigl(a(x_i)u_i^2-b(x_i)u_i\bigr)\frac{\partial_{u_i}g(x_i,u_
i)}{g(x_i,u_i)},\nonumber
\end{eqnarray}
where
$g\dvtx E\times[0,\infty)\rightarrow(0,1]$ has the property
that there exists $r_g$ such that $g(z,v)=1$ for $v>r_g$,
and $f(x,u)=\prod_ig(x_i,u_i)$. We can identify the state space
of the corresponding process with a subset of
$(E\times[0,\infty))^{\infty}$ or, since order is not important, with a
subset of $\mathcal{N}(E\times[0,\infty))$, the counting measures on
$E\times[0,\infty)$. Define
\[
\mathcal{N}_f\bigl(E\times[0,\infty)\bigr)=\bigl\{\mu\in\mathcal{N}\bigl(E\times
[0,\infty
)\bigr)\dvtx \mu\bigl(E\times[0,u]\bigr)<\infty\ \forall 0<u<\infty\bigr\},
\]
where the topology for $\mathcal{N}_f(E\times\infty)$ is given by the
requirement that $\mu_n\rightarrow\mu$ if and only if $\int f\,d\mu_
n\rightarrow\int f\,d\mu$ for
all $f\in\overline{C}(E\times[0,\infty)$ for which there exists $u_
f>0$ such that
$f(x,u)=0$ for $u\geq u_f$. (See Section \ref{convres} for a
discussion of the appropriate topology to use in the
infinite measure setting.)

As $r\rightarrow\infty$, the particle process converges to a process in
which particle locations evolve as independent Markov
processes with generator $B$, levels satisfy~(\ref{leveq}), a
particle with level $U_i(t)$
gives birth to new particles at its location $X_i(t)$ and
level in the interval $[U_i(t)+c,U_i(t)+d]$ at rate
$2a(X_i(t))(d-c)$. A~particle dies when its level hits $\infty$.
The level of a particle born at time $t_0$ (or in
the initial population, if $t_0=0$) with
initial level $U_i(t_0)$ satisfies
\[
U_i(t)=\frac{U_i(t_0)e^{-\int_{t_0}^tb(X_i(s))\,ds}}{1-U_i(t_0)\int_{
t_0}^te^{-\int_{t_0}^vb(X_i(s))\,ds}a(X_i(v))\,dv},
\]
until it hits infinity. If $b\leq0$ and $a$ is bounded away
from zero, then $U_i$ will hit infinity in finite time.

If we extend the path $X_i$ back along its
ancestral path to time zero, we would have
\[
U_i(t)\geq\frac{u_0e^{-\int_0^tb(X_i(s))\,ds}}{1-u_0\int_0^te^{-
\int_0^vb(X_i(s))\,ds}a(X_i(v))\,dv},
\]
where $u_0$ is the level of the particle's ancestor at time
zero, and $X_i(s)$ is the position of the ancestor at time
$s\leq t$. Since we are assuming that the initial position of an
offspring is that of the parent, $X_i$ is a solution of the
martingale problem for $B$.
\begin{proposition}\label{bprest}
Let $(X,U)$ be a solution of the martingale problem for $A$
given by (\ref{eq33}). Let $(X^r,U^r)$ consist of the subset of
particles for which $U_i<r$, that is,
\[
\sum_{}\delta_{(X_j^r(t),U_j^r(t))}=\sum\delta_{(X_i(t),U_i(t))}
{\mathbf1}_{[0,r)}(U_i(t)).
\]
If $a(z)r-b(z)\geq0$ for all $z\in E$, then $(X^r,U^r)$ is a solution
of the martingale problem for $A_r$ given by (\ref{fingen}).\vadjust{\goodbreak}
\end{proposition}
\begin{remark}
The condition $a(z)r-b(z)\geq0$ for all $z\in E$ ensures that any particle
that is above level $r$ at time $t$ will stay above level
$r$ at all future times.
\end{remark}
\begin{pf*}{Proof of Proposition \ref{bprest}}
The proposition follows by the observation that $A_r$ can
be obtained from $A$ by restricting the domain to
$f(x,u)=\prod_ig(x_i,u_i)$ for which $g(z,v)=1$ for $v\geq r$.
\end{pf*}

%s3.3 ###
\subsection{The genealogy}\label{geneal2}
Assume for the moment that $a$ and $b$ do not depend on $x$.
If the location process has the property that at the time
of a birth only the offspring has the same location as
the parent (e.g., if the location process is
Brownian motion), then the full genealogy can be
recovered from knowledge of the levels and locations.
The collection of ancestors at
time $t<T$ of the
particles alive at time $T$ is $\{(X_i(t),U_i(t))\dvtx U_i(t)<u_T(t)\}$,
where $u_T(t)$ is given by (\ref{popbnd}). The number of
particles in this collection is nondecreasing in $t$ and
increases only by jumps of $+1$. The parent of the new
particle is identifiable by the fact that only the parent
and the offspring will be at the same location.

If $a$ and $b$ depend on $x$, the full genealogy is still
determined by the locations and levels of the particles,
but recovering the genealogy is more complicated since
it may not be possible to tell whether or not a particle
$(X_i(t),U_i(t))$ has descendants alive at time $T>t$ just
from information available at time $t$. However, some
easy observations can be made. For example, if
$\inf_xa(x)>0$ and $\sup_xb(x)<\infty$, then for $t<T$, all particles
alive at
time $T$ are descendants of finitely many particles alive
at time $t$.

%s3.4 ###
\subsection{The measure-valued limit}\label{sec-mes-val}
The generator for a Dawson--Watanabe superprocess is
typically of the form
\[
C\widehat{f}(\mu)=\exp\{-\langle h,\mu\rangle\}\int_E\bigl(-Bh(y)-F(h
(y),y)\bigr)\mu(dy),
\]
for $\widehat{f}(\mu)=\exp\{-\langle h,\mu\rangle\}$, where $\langle
h,\mu\rangle=\int_Eh(y)\mu(dy)$ and $h$ is
an appropriate function in $\mathcal{D}(B)$ (see, e.g.,
Theorem 9.4.3 of \cite{EK86}). For superprocesses
arising from branching models with offspring
distributions having finite variances, $F$ should be of the
form $F(h(y),y)=-a(y)h(y)^2+b(y)h(y)$.

For $\mu\in\mathcal{M}_f(E)$, let $\alpha(\mu,dx\times du)$ be the
distribution of a
Poisson random measure on $E\times[0,\infty)$ with mean measure
$\mu\times\varLambda$. Then setting $h(y)=\int_0^{\infty}(1-g(y,v))\,dv$,
\begin{eqnarray*}
\widehat{f}(\mu)&=&\alpha f(\mu)=\int f(x,u)\alpha(\mu,dx\times
du)\\
&=&\exp\biggl\{\int_E\int_0^{\infty}\bigl(g(y,v)-1\bigr)\,dv\,\mu(dy)\biggr\}\\
&=&\exp\{-\langle
h,\mu\rangle\}.
\end{eqnarray*}
Using Lemma
\ref{lemma2}, we have
\begin{eqnarray*}
\alpha Af(\mu)&=&\int_E\int_0^{\infty}Bg(y,v)\,dv\,\mu(dy)\exp\{-\langle
h,\mu\rangle\}\\[-0.53pt]
&&{} +\int_E\int_0^{\infty}2a(y)g(y,z)\int_z^{\infty}\bigl(g(y,v)-
1\bigr)\,dv\,dz\,\mu(dy)\exp\{-\langle h,\mu\rangle\}\\[-0.53pt]
&&{} +\int_E\int_0^{\infty}\bigl(a(y)v^2-b(y)v\bigr)\,\partial_vg(y,v)\,dv\,\mu
(dy)\exp\{-\langle h,\mu\rangle\}\\[-0.53pt]
&=&\int_E\int_0^{\infty}Bg(y,v)\,dv\,\mu(dy)\exp\{-\langle h,\mu\rangle
\}\\[-0.53pt]
&&{} +\int_E\int_0^{\infty}2a(y)g(y,z)\int_z^{\infty}\bigl(g(y,v)-
1\bigr)\,dv\,dz\,\mu(dy)\exp\{-\langle h,\mu\rangle\}\\[-0.53pt]
&&{} -\int_E\int_0^{\infty}\bigl(2a(y)v-b(y)\bigr)\bigl(g(y,v)-1\bigr)\,dv\,\mu(dy)\exp
\{-\langle h,\mu\rangle\}\\[-0.53pt]
&=&\int_E\int_0^{\infty}Bg(y,v)\,dv\,\mu(dy)\exp\{-\langle h,\mu\rangle
\}\\[-0.53pt]
&&{} +\int_E\int_0^{\infty}2a(y)g(y,z)\int_z^{\infty}\bigl(g(y,v)-
1\bigr)\,dv\,dz\,\mu(dy)\exp\{-\langle h,\mu\rangle\}\\[-0.53pt]
&&{} -\int_E\int_0^{\infty}2a(y)\int_z^{\infty}\bigl(g(y,v)-1\bigr)\,dv\,\mu
(dy)\exp\{-\langle h,\mu\rangle\}\\[-0.53pt]
&&{} +\int_E\int_0^{\infty}b(y)\bigl(g(y,v)-1\bigr)\,dv\,\mu(dy)\exp\{-\langle
h,\mu\rangle\}\\[-0.53pt]
&=&\int_E\int_0^{\infty}Bg(y,v)\,dv\,\mu(dy)\exp\{-\langle h,\mu\rangle
\}\\[-0.53pt]
&&{} +\int_Ea(y)\biggl(\int_0^{\infty}\bigl(g(y,v)-1\bigr)\,dv\biggr)^2\mu
(dy)\exp\{-\langle h,\mu\rangle\}\\[-0.53pt]
&&{} +\int_E\int_0^{\infty}b(y)\bigl(g(y,v)-1\bigr)\,dv\,\mu(dy)\exp\{-\langle
h,\mu\rangle\}\\[-0.53pt]
&=&\int_E\bigl(-Bh(y)+a(y)h(y)^2-b(y)h(y)\bigr)\mu(dy)\exp\{-\langle
h,\mu\rangle\}
=C\widehat{f}(\mu).
\end{eqnarray*}
But $C$ is the generator for a superprocess, so for each
$\mu\in\mathcal{M}_f(E)$, there exists a solution $Z$ of the martingale
problem for $C$ with $Z(0)=\mu$ and hence a solution $(X,U)$
of the martingale problem for $A$ with initial distribution
$\alpha(\mu,\cdot)$.

The mapping $\gamma\dvtx \mathcal{N}_f(E\times[0,\infty))\rightarrow
\mathcal{M}_
f(E)$
used in the application of Theorem \ref{mf} is given by
\[
\gamma\Bigl(\sum\delta_{(x_i,u_i)}\Bigr)=\cases{
\displaystyle \lim_{r\rightarrow\infty}\frac1r\sum_{u_i\leq r}\delta_{x_i},&\quad
if the measures converge\vspace*{-10pt}\cr
&\quad in the weak topology,\vspace*{2pt}\cr
\mu_0, &\quad otherwise,}
\]
where $\mu_0$ is a fixed element of $\mathcal{M}_f(E)$.
The solution $\sum\delta_{(X_i(t),U_i(t))}$ of the martingale problem for
$A$ is a conditionally Poisson random measure (see
Section \ref{sectcp}) with Cox
measure $Z(t)$. Consequently, the particles determine $Z$ by
%
%e3.6 ###
\begin{equation}\label{bplim}
Z(t)=\lim_{r\rightarrow\infty}\frac1r\sum_{U_i(t
)\leq r}\delta_{X_i(t)}.
\end{equation}

Since by Proposition \ref{bprest} and Theorem
\ref{finbmp}, $Z^r(t)=\frac1r\sum_{U_i(t)\leq r}\delta_{X_i(t)}$ is
the normalized
empirical measure for a branching Markov process,
(\ref{bplim}) gives the convergence of the normalized
branching Markov process to the corresponding
Dawson--Watanabe superprocess (cf. \cite{Wat68,Daw75}).

%s4 ###
\section{Examples and extensions}\label{secteg}

%s4.1 ###
\subsection{A model with immigration} The simplest
immigration process assumes that the space--time point
process giving the arrival times and locations of the
immigrants is a Poisson process. Assuming temporal
homogeneity, immigration is
introduced by adding the generator of a
space--time-level Poisson random measure. Let $\nu$ be
the intensity of immigration, that is, $\nu(A)\Delta t$ is
approximately the
probability that an individual immigrates into $A\subset E$ in
a time interval of length~$\Delta t$. The generator becomes
%
%e4.1 ###
\begin{eqnarray}\label{genim}
Af(x,u)&=&f(x,u)\sum_i\frac{Bg(x_i,u_i)}{g(x_i,u_i)}\nonumber\\
&&{}+f(x,u)\sum_
i2a(x_i)\int_{u_i}^{r_g}\bigl(g(x_i,v)-1\bigr)\,dv\nonumber\\[-8pt]\\[-8pt]
&&{}+f(x,u)\int_0^{r_g}\int_E\bigl(g(z,v)-1\bigr)\nu(dz)\,dv\nonumber\\
&&{}+f(x,u)\sum_i\bigl(a(x_i)u_i^2-b(x_i)u_i\bigr)\,\frac{\partial_{
u_i}g(x_i,u_i)}{g(x_i,u_i)}.\nonumber
\end{eqnarray}
Noting that the generator for finite $r$ is obtained from $A$
by restricting the domain to the collection of $g$ with
$r_g\leq r$, if $ra(x)-b(x)\geq0$, for all $x$, the generator of
the corresponding branching Markov process with
immigration is
\begin{eqnarray*}
C_r\widehat{f}(x,n)&=&\sum_{i=1}^nB_{x_i}\widehat{f}(x,n)+\sum_{i=1}^nr
a(x_i)\bigl(\widehat{f}\bigl(b(x|x_i),n+1\bigr)-\widehat{f}(x,n)\bigr)\\
&&{} +\sum_{i=1}^n\bigl(ra(x_i)-b(x_i)\bigr)\bigl(\widehat{f}\bigl(d(x|x_i),n-1\bigr)-\widehat{
f}(x,n)\bigr)\\
&&{} +\int_E\bigl(\widehat{f}\bigl(b(x|z),n+1\bigr)-\widehat{f}(x,n)\bigr)\nu(dz).
\end{eqnarray*}

For $r=\infty$, setting $h(x)=\int_0^{\infty}(1-g(x,v))\,dv$ as before, the
generator for the measure-valued process is
\[
C\widehat{f}(\mu)=\alpha Af(\mu)=(\langle-Bh+ah^2-bh,\mu\rangle
-\langle h,\nu\rangle)\exp\{-\langle h,\mu\rangle\},
\]
for $\widehat{f}(\mu)=\exp\{-\langle h,\mu\rangle\}$.

Early results on branching Markov processes with
immigration include \cite{Her78,Kul79}. Work on
limiting
measure-valued processes with immigration includes
\cite{EG93,LiZ92,LLW93,LiWa99,Sta03}.

%s4.2 ###
\subsection{Conditioning on nonextinction}
In the limiting model considered in Section \ref{rlim}, let $a$
and $b$ be constant and
$b<0$. Let $\tau$ be the time of
extinction and let $U_{*}(0)$ be the minimum of the initial
levels. Then $\tau$ is the solution of $1-U_{*}(0)\frac ab(1-e^{
-b\tau})=0$,
so
\[
\tau=-\frac1b\log\frac{U_{*}(0)a-b}{U_{*}(0)a}.
\]
If $Z(0)=\mu_0$, then $U_{*}(0)$ is exponentially distributed with
parameter $\mu_0(E)$ and
\[
P\{\tau>T\}=P\biggl\{U_{*}(0)<\frac b{a(1-e^{-bT})}\biggr\}=1-\exp\biggl\{-\frac{
b\mu_0(E)}{a(1-e^{-bT})}\biggr\}.
\]
The case $b=0$ is obtained by passing to the limit so
that $P\{\tau>T\}=1-\exp\{-\mu_0(E)/(aT)\}$.

As in Section \ref{cndne}, conditioning on
$\{\tau>T\}$ and letting $T\rightarrow\infty$ is equivalent to
conditioning on
the initial Poisson process having a level at zero. The
resulting generator becomes
%
%e4.2 ###
\begin{eqnarray}\label{gencnd}
Af(x,u)&=&f(x,u)\sum_i\frac{Bg(x_i,u_i)}{g(x_i,u_i)}\nonumber\\
&&{}+f(x,u)\sum_{
u_i>0}2a\int_{u_i}^{r_g}\bigl(g(x_i,v)-1\bigr)\,dv\nonumber\\[-8pt]\\[-8pt]
&&{} +f(x,u)2a\int_0^{r_g}\bigl(g(x_0,v)-1\bigr)\,dv\nonumber\\
&&{} +f(x,u)\sum_{u_i>0}(au_i^2-bu_i)\,\frac{\partial_{u_i}g(x_
i,u_i)}{g(x_i,u_i)},\nonumber
\end{eqnarray}
where the $u_i$ are the nonzero levels,
and the generator for the conditioned measure-valued process is
given by setting
\begin{eqnarray*}
\alpha_0f(\mu)&=&\int f(x,u)\alpha_0(\mu,dx\times du)\\
&=&\frac1{|
\mu|}\int_Eg(z,0)\mu(dz)\exp\{-\langle h,\mu\rangle\}
\end{eqnarray*}
and
\begin{eqnarray*}
\alpha_0Af(\mu)&=&\langle-Bh(y)+ah(y)^2-b(y)h(y),\mu\rangle\\
&&{}\times\frac
1{|\mu|}\int_Eg(z,0)\mu(dz)\exp\{-\langle h,\mu\rangle\}\\
&&{}+\frac1{|\mu|}\int_E\bigl(Bg(z,0)-2ag(z,0)h(z)\bigr)\mu(dz)\exp\{-\langle
h,\mu\rangle\}.
\end{eqnarray*}
Note that $\alpha_0$ is the distribution of $\xi=\delta_{0,Z_0}+
\sum_{i=1}^{\infty}\delta_{(U_i,Z_i)}$,
where $\{U_i,i\geq1\}$ is a Poisson process with intensity $\mu
(E)$, and
$Z_0,Z_1,\ldots$ are i.i.d. with distribution $\mu(\cdot)/\mu(E)$.

For earlier work, see \cite{Ev93,EP90,HH81,LN68,Mel82}.
In particular, the particle at level zero in the
construction above is the ``immortal
particle'' of Evans \cite{Ev93}.

%s4.3 ###
\subsection{Conditioning on extinction}
As in Section \ref{cnde}, if $a$ and $b$ are constant and $b>0$, then
conditioning on
extinction is equivalent to conditioning on $U_{*}(0)>\frac ba$.
Defining $V_i(t)=U_i(t)-\frac ba$, the generator for the
conditioned process is
%
%e4.3 ###
\begin{eqnarray}\label{gencnd3}
Af(x,v)&=&f(x,v)\sum_i\frac{Bg(x_i,v_i)}{g(x_i,v_i)}\nonumber\\
&&{}+f(x,v)\sum_{i
>0}2a\int_{v_i}^{r_g}\bigl(g(x_i,v)-1\bigr)\,dv\\
&&{} +f(x,v)\sum_{i>0}(av_i^2-bv_i)\,\frac{\partial_{v_i}g(x_i,v_i)}{g(
x_i,v_i)},\nonumber
\end{eqnarray}
and the generator of the measure-valued process is
\[
C\widehat{f}(\mu)=\int_E\bigl(-Bh(y)+ah(y)^2+bh(y)\bigr)\mu(dy)\exp
\{-\langle h,\mu\rangle\},
\]
for $\widehat{f}(\mu)=\exp\{-\langle h,\mu\rangle\}$. In other words,
conditioning a
supercritical process on extinction replaces the
supercritical process by a subcritical one. This result
is originally due to Evans and O'Connell \cite{EvOC94}.

%s4.4 ###
\subsection{Models with multiple simultaneous births}
We now consider continuous-time, branching Markov
processes with general offspring distributions.
The general theory of branching Markov processes was
developed by Ikeda, Nagasawa and Watanabe in a long
series of papers \cite{INW65,INW68a,INW68b,INW69}
following earlier work by several authors. The
particle representation is substantially more complicated
and passage to the infinite population limit more
delicate.\vadjust{\goodbreak}

As above, the particles move independently in $E$ according to
a generator $B$. A particle at position $x\in E$ with level $u$
gives birth to $k$ offspring at rate
$(k+1)a^{(r)}_k(x)(r-u)^kr^{-(k-1)}$. New particles have the
location of the parent,
but their levels are
uniformly distributed on $[u,r)$.
Then for
$f(x,u,n)=\prod_{i=1}^ng(x_i,u_i)$,
%
%e4.4 ###
\begin{eqnarray}\label{multgen}\qquad
&&A_rf(x,u,n)\nonumber\\
&&\qquad =f(x,u,n) \sum_{i=1}^n\frac{Bg(x_i,u_i)}{g(x_i,u_i)}
\nonumber\\
&&\qquad\quad{} +f(x,u,n) \sum_{i=1}^n \sum_{k=1}^{\infty}\frac{(k+1)
a_k^{(r)}(x_i)}{r^{k-1}} \nonumber\\
&&\qquad\quad\hspace*{89.1pt}{}\times\int_{[u_i,r)^k}\Biggl[\Biggl(\prod_{l=1}^
kg(x_i,v_l)\Biggr)-1\Biggr] \,dv_1\cdots dv_k\\
&&\qquad\quad{} +f(x,u,n)\sum_{i=1}^n \Biggl(\sum_{k=1}^{\infty}r^2a_k^{
(r)}(x_i)\biggl[\biggl(1-\frac{u_i}r\biggr)^{k+1}\nonumber\\
&&\qquad\quad\hspace*{144.2pt}{}-1+(k+1)\frac{u_
i}r\biggr]\nonumber\\
&&\qquad\quad\hspace*{180pt}{}-b(x_i)u_i\Biggr)\,\frac{\partial_{u_i}g(x_i,u_i)}{g(x_i
,u_i)}.\nonumber
\end{eqnarray}
Now, the levels satisfy the equation
%
%e4.5 ###
\begin{eqnarray}\label{levdif}
\dot{U}_i(t)&=&\sum_{k=1}^{\infty}r^2a_k^{(r)}(X_i
(t))\biggl[\biggl(1-\frac{U_i(t)}r\biggr)^{k+1}-1+(k+1)\frac{U_i(
t)}r\biggr]\nonumber\\[-8pt]\\[-8pt]
&&{}-b(X_i(t))U_i(t).\nonumber
\end{eqnarray}

Defining $g(x)=\frac1r\int_0^rg(x,v) \,dv$ and integrating
(\ref{multgen}) with
respect to $\alpha(n,du)$, the uniform measure on
$[0,r]^n$, we have that
%
%e4.6 ###
\begin{eqnarray}\label{mgenprj}
&&\int A_rf(x,u,n) \alpha(n,du)\nonumber\\
&&\qquad =C_rf(x,n)\nonumber\\
&&\qquad =f(x,n) \sum_{i=1}^n\frac{Bg(x_i)}{g(x_i)}\\
&&\qquad\quad{}+f(x,n) \sum_{
i=1}^n \sum_{k=1}^{\infty} ra_k^{(r)}(x_i) [g(x_i)^k-1
]\nonumber\\
&&\qquad\quad{} +f(x,n) \sum_{i=1}^n \Biggl(r\sum_{k=
1}^{\infty} k a_k^{(r)}(x_i)-b(x_i)\Biggr)\biggl[\frac1{g(x_i)}
-1\biggr],\nonumber
\end{eqnarray}
which is the generator of a branching
process with multiple births with birth rates
$ra_k^{(r)}(\cdot)$, death rate $r\sum_{k=1}^{\infty} k a_k^{(
r)}(\cdot)-b(\cdot)$ (provided this
expression is nonnegative) and particles
moving according to the generator $B$. The analog of
Theorem \ref{finbmp} holds with $a(x_i)$ replaced by
$\sum_k(k+1)a^{(r)}_k(x_i)$ in the definition of $\psi$ and $\widetilde{
\psi}$.

We define
\[
\Lambda^{(r)}(x,u)=\sum_{k=1}^{\infty}r(k+1)a_k^{(r)}(x)\biggl[1
-\biggl(1-\frac ur\biggr)^k\biggr]
\]
and assume that
%
%e4.7 ###
\begin{eqnarray}\label{lamdef}
&&\lim_{r\rightarrow\infty}\Lambda^{(r)}(x,u)\nonumber\\
&&\qquad =\lim_{r\rightarrow\infty}\sum_{k=1}^{\infty}r(k+1)a_k^{
(r)}(x)\Biggl[\sum_{l=1}^k\pmatrix{k\cr l}(-1)^{l+1}\biggl(\frac ur
\biggr)^l\Biggr]\\
&&\qquad \equiv\Lambda(x,u)\nonumber
\end{eqnarray}
exists uniformly for $x\in E$ and $u$ in bounded intervals. This
condition is essentially (9.4.36) of \cite{EK86}.

Observe that
\[
\int_0^u\Lambda^{(r)}(x,v)\,dv=\sum_{k=1}^{\infty}r^2a_k^{(r)}(x)
\biggl[\biggl(1-\frac ur\biggr)^{k+1}-1+(k+1)\frac ur\biggr],
\]
so that (\ref{levdif}) becomes
%
%e4.8 ###
\begin{equation}\label{levdif2}
\dot{U}_i(t)=\int_0^{U_i(t)}\Lambda^{(r)}(X_i(s)
,v)\,dv-b(X_i(t))U_i(t)
\end{equation}
and
\[
\int_0^r\Lambda^{(r)}(x,v)\,dv=\sum_{k=1}^{\infty}r^2ka_k^{(r)}(x
),
\]
so that the death rate for the branching process can be
written as $r^{-1}\int_0^r\Lambda^{(r)}(x,\break v)\,dv-b(x)$.

As in \cite{K00},
%
%e4.9 ###
\begin{eqnarray}\label{subint}
\partial^m\Lambda^{(r)}(x,u)&\equiv&\frac{\partial^m}{\partial u^
m}\Lambda^{(r)}(x,u)\nonumber\\
&=&(-1)^{m+1}\sum_{k=m}^{\infty}a_k^{(r)}(x)\frac{(k+1)k\cdots(
k-m+1)}{r^{m-1}}\\
&&\hspace*{57.6pt}{}\times\biggl(1-\frac ur\biggr)^{k-m}.\nonumber
\end{eqnarray}
Each of the derivatives is a monotone function of $u$,
and since
\begin{eqnarray*}
&&\partial^{m-1}\Lambda^{(r)}(x,b)-\partial^{m-1}\Lambda^{(r)}(x,a)\\
&&\qquad=\int_a^b\partial^m\Lambda^{(r)}(x,u)\,du,
\end{eqnarray*}
it follows by the convergence of $\Lambda^{(r)}$ and
induction on $m$ that each $\partial^m\Lambda^{(r)}(x,u)$ converges
uniformly in $u$ on bounded intervals that are bounded
away from zero. Consequently, $\Lambda$ is infinitely differentiable in $
u$,
and for $0<u_1<{u_2<\infty}$,
\[
\lim_{r\rightarrow\infty}\sup_{u_1 \leq u\leq u_2}\bigl|\partial^m\Lambda^{
(r)}(x,u)-\partial^m\Lambda(x,u)\bigr|=0.
\]
The fact that the derivatives alternate in sign implies
that $\partial^1\Lambda(x,\cdot)$ is completely monotone and hence
can be
represented as
\[
\partial^1\Lambda(x,u)=\int_0^{\infty}e^{-uz}\widehat{\nu}(x,dz)
\]
for some $\sigma$-finite measure $\widehat{\nu}(x,\cdot)$.
Writing $\widehat{\nu}(x,\cdot)=2a_0(x)\delta_0+\nu(x,\cdot)$ with $
\nu(x,\{0\})=0$,
\[
\Lambda(x,u)=2a_0(x)u+\int_0^{\infty}z^{-1}(1-e^{-uz})\nu(x,d
z).
\]

Let $g$ satisfy $g(x,v)=1$, for
$v\geq u_g$, and define
\[
h(x,u)=\int_u^{u_g}\bigl(1-g(x,v) \bigr)\,dv.
\]
If $r>u_g$ and there are $n$ particles below level $r$, then
(\ref{multgen}) may be written as
\begin{eqnarray*}
&&A_rf(x,u,n)\\
&&\qquad =f(x,u,n) \sum_{i=1}^n\frac{Bg(x_i,u_i)}{g(x_i,u_i)}\\
&&\qquad\quad{} +f(x,u,n)\\
&&\qquad\quad\hspace*{10.4pt}{}\times\sum_{i=1}^n\sum_{k=1}^{\infty}r (k+1)\\
&&\qquad\quad\hspace*{52.83pt}{}\times a_k^{(r)}
(x_i)\biggl\{\biggl(1-\frac{u_i+h(x_i,u_i)}r\biggr)^k-\biggl(1-\frac{
u_i}r\biggr)^k\biggr\}\\
&&\qquad\quad{} +f(x,u,n)\sum_{i=1}^n \biggl(\int_0^{u_i}\Lambda^{(r)}(x_
i,v)\,dv-b(x_i)u_i\biggr)\,\frac{\partial_{u_i}g(x_i,u_i)}{g(x_i,u_i
)}.
\end{eqnarray*}
Then, by (\ref{lamdef}) and the definition of $h$, we have
\begin{eqnarray*}
&&Af(x,u)\\
&&\qquad\equiv\lim_{r\rightarrow\infty}A_rf(x,u)\\
&&\qquad=f(x,u) \sum_i\frac{Bg(x_i,u_i)}{g(x_i,u_i)}\\
&&\qquad\quad{} +f(x,u)\sum_i\bigl(\Lambda(x_i,u_i)-\Lambda\bigl(x_i,u_i+h(x_i,u_
i)\bigr)\bigr)\\
&&\qquad\quad{} +f(x,u)\sum_i\biggl(\int_0^{u_i}\Lambda(x_i,v)\,dv-b(x_i)u_
i\biggr)\,\frac{\partial_{u_i}g(x_i,u_i)}{g(x_i,u_i)}\\
&&\qquad=f(x,u) \sum_i\frac{Bg(x_i,u_i)}{g(x_i,u_i)}\\
&&\qquad\quad{} +f(x,u)\sum_i2a_0(x_i)\int_{u_i}^{\infty}\bigl(g(x_i,v)-1\bigr)\,dv\\
&&\qquad\quad{} +f(x,u)\sum_i\int_0^{\infty}\bigl(e^{z\int_{u_i}^{\infty}(g(x_
i,v)-1)\,dv}-1\bigr)z^{-1}e^{-zu_i}\nu(x_i,dz)\\
&&\qquad\quad{} +f(x,u)\sum_i\biggl(a_0(x_i)u_i^2-b(x_i)u_i\\
&&\qquad\quad\hspace*{66.3pt}{} +\int_0^{\infty}z^{-1}\bigl(u_i-z^{-1}(1-e^{
-u_iz})\bigr)\nu(x_i,dz)\biggr)\,\frac{\partial_{u_i}g(x_i,u_i)}{g(x_i,
u_i)}.
\end{eqnarray*}
Note that the second term on the right-hand side has the same
interpretation as the second term on the right-hand side of
(\ref{eq33}). To understand the third term,
recall that if $\xi=\sum_i\delta_{\tau_i}$ is a Poisson process on $
[0,\infty)$ with
parameter $\lambda$, then
\[
E\Bigl[\prod g(\tau_i)\Bigr]=e^{\lambda\int_0^{\infty}(g(v)-1)\,dv}.
\]
Consequently, the third term determines bursts of
simultaneous offspring
at the location $x_i$ of the parent and with levels forming a
Poisson process with intensity $z$ on $[u_i,\infty)$.

Setting $h(x)=h(x,0)=\int_0^{\infty}(1-g(x,v))\,dv$ and
$\widehat{f}(\mu)=\exp\{-\langle h,\mu\rangle\}$,
\begin{eqnarray*}
C\widehat{f}(\mu)&=&\alpha Af(\mu)\\
&=&\int_E\biggl(-Bh(y)+\int_0^{h(y)}\Lambda(y,z)\,dz-b(y)h(y)
\biggr)\mu(dy)\exp\{-\langle h,\mu\rangle\}.
\end{eqnarray*}

Based on the above calculations, we have the following
theorem.
\begin{theorem}
Let $B\subset\overline{C}(E)\times\overline{C}(E)$ satisfy Condition \ref
{bmpc}, and let
the martingale problem for $B$ be well posed.
Assume that for $r\geq r_0$,
\[
\inf_x\biggl(r\sum_{k=1}^{\infty} k a_k^{(r)}(x)-b(x)\biggr)\geq
0
\]
and that the convergence in (\ref{lamdef}) is uniform in $x$.
Let $K(0)$ be a finite random measure on $E$, and let $\xi^r$ be a
solution of the martingale problem for $A_r$ such that $\xi^r(0)$
is conditionally Poisson on $E\times[0,r]$ with mean measure
$K(0)\times\varLambda$. Then $\xi^r\Rightarrow\xi$ where $\xi$ is a
solution of the
martingale problem for $A$.
\end{theorem}
\begin{pf}
For $r>q$, let $\xi^{(q),r}$ denote the restriction of $\xi^r$ to $
E\times[0,q]$ and
similarly for $\xi^{(q)}$. It is enough to prove that
$\xi^{(q),r}\Rightarrow\xi^{(q)}$ for each $q>r_0$. The generator for $
\xi^{(q),r}$ is
the restriction of $A_r$ to functions $f\in\mathcal{D}(A_r)$ such that the
corresponding $g$ satisfies $g(x,u)=1$ for $u\geq q$. For $f$ of
this form, by (\ref{subint}), $A_{r,q}f=A_rf$ satisfies
\begin{eqnarray*}
&&A_{r,q}f(x,u,n)\\
&&\qquad =f(x,u,n) \sum_{i=1}^n\frac{Bg(x_i,u_i)}{g(x_i,u_i)}
\\
&&\qquad\quad{} +f(x,u,n) \sum_{i=1}^n \sum_{m=1}^{\infty}\frac1{m!}(
-1)^{m+1}\partial^m\Lambda^{(r)}(x_i,q)\\
&&\qquad\quad\hspace*{91.1pt}{}\times\int_{[u_i,q)^m}\Biggl[
\Biggl(\prod_{l=1}^mg(x_i,v_l)\Biggr)-1\Biggr] \,dv_1\cdots dv_m\\
&&\qquad\quad{} +f(x,u,n)\sum_{i=1}^n \biggl(\int_0^{u_i}\Lambda^{(r)}(x
,v)\,dv-b(x_i)u_i\biggr)\,\frac{\partial_{u_i}g(x_i,u_i)}{g(x_i,u_i)},
\end{eqnarray*}
and
the corresponding branching Markov process has
generator
\begin{eqnarray*}
&&C_{r,q}f(x,n)\\
&&\qquad =f(x,n) \sum_{i=1}^n\frac{Bg(x_i)}{g(x_i)}\\
&&\qquad\quad{} + f(x,n) \sum_{
i=1}^n \sum_{m=1}^{\infty}\frac1{(m+1)!}(-1)^{m+1}\partial^m\Lambda^{
(r)}(x_i,q)q^m[g(x_i)^m-1]\\
&&\qquad\quad{} +f(x,n) \sum_{i=1}^n \biggl(\frac1q\int_
0^q\Lambda^{(r)}(x_i,v)\,dv-b(x_i)\biggr)\biggl[\frac1{g(x_i)}-1
\biggr].
\end{eqnarray*}
The convergence of $\xi^{(q),r}$ follows by the convergence
assumptions on $\Lambda^{(r)}$.
\end{pf}

The measure $\nu(x,\cdot)$ is nonzero only if the offspring
distribution has a ``heavy tail.'' If
$a_k^{(r)}(x)=a_k(x)$
and
\[
\sum_{k=1}^{\infty}(k+1) k a_k(x)<\infty,
\]
then
\[
\Lambda(x,u)=\lim_{r\rightarrow\infty}\sum_{k=1}^{\infty}r(k+1
)a_k(x)\biggl[1-\biggl(1-\frac ur\biggr)^k\biggr]=\sum_{k=1}^{\infty}
(k+1)ka_k(x)u
\]
and
\begin{eqnarray*}
Af(x,u)&=&f(x,u)\sum_{i }\frac{Bg(x_i,u_i)}{g(x_i,u_i)}\\
&&{} +f(x,u)\sum_{i }\sum_{k=1}^{\infty}(k+1)ka_k(x_i)\int_{
u_i}^{u_g}[g(x_i,v)-1] \,dv\\
&&{} +f(x,u)\sum_{i }\Biggl(\sum_{k=1}^{\infty}\frac{(k+1)
k a_k(x_i)}2 u_i^2-b(x_i)u_i\Biggr)\, \frac{\partial_{u_i}g(x_i
,u_i)}{g(x_i,u_i)},
\end{eqnarray*}
which is essentially (\ref{eq33}).

For scalar branching processes with general offspring
distributions, convergence to possibly discontinuous
continuous state branching processes was proved by
Grimvall \cite{Gri74} (see \cite{EK86}, Section 9.1).
Convergence in the measure-valued setting is given in
\cite{Wat68} and \cite{Daw75} for offspring distributions with finite
second moment and more generally in
\cite{EK86}, Theorem 9.4.3. Fitzsimmons \cite{Fit88}
gives a very general construction of these processes.

If $\nu$ is not zero, then the genealogy of the process is
much more complicated than that described in Sections
\ref{geneal} and \ref{geneal2}. Assume that $\Lambda$ and $b$ do
not depend on~$x$, and define
\[
\widehat{\Lambda}(u)\equiv\int_0^u\Lambda(v)\,dv-bu=2a_0u^2-bu+\int_
0^{\infty}z^{-2}(zu-1+e^{-uz})\nu(dz).
\]
If $u_T$ satisfies
\[
\dot{u}_T(t)=\widehat{\Lambda}(u_T(t)),
\]
$u_T(t)<\infty$ for $t<T$ and $\lim_{t\rightarrow T-}u_T(t)=\infty$,
then it is
still the case that the collection of ancestors at
time $t<T$ of the population alive at time $T$ is
$\{(X_i(t),U_i(t))\dvtx U_i(t)<u_T(t)\}$, but $u_T$ may not exist. In fact,
since if $\dot{u}=\widehat{\Lambda}(u)$
\[
\int_{u(t_1)}^{u(t_2)}\frac1{\widehat{\Lambda}(v)}\,dv=t_2-t_1,
\]
$u_T$ exists if and only if
\[
\int_u^{\infty}\frac1{\widehat{\Lambda}(v)}\,dv<\infty
\]
for $u$ sufficiently large, which always holds if $a_0>0$.
In the critical and subcritical
cases, this condition is equivalent to extinction with
probability one as was noted by Bertoin and Le Gall
(\cite{BL06}, page 167).

This
\textit{finite ancestry} property or \textit{coming down from infinity }
of the genealogy has been studied for a variety of
population models. See \cite{Sch00b} and \cite{BBL09} for
results in the Fleming--Viot setting. The equivalence of
the conditions for Fleming--Viot and Dawson--Watanabe
processes is given in (\cite{BL06}, page 171).

The argument in Section \ref{harthm} can undoubtedly be
extended to the present setting. This development will
be carried out elsewhere.

%s4.5 ###
\subsection{Model with exponentially distributed levels}
The discrete models that we have considered have been
formulated with levels that are uniformly distributed on
an interval. That is not necessary, and other
distributions may be convenient in other contexts. We
illustrate this flexibility by
formulating a model for a simple branching process with
levels that are exponentially distributed. The dynamics
of the levels change, and the correct dynamics are
determined by essentially working backwards from the
answer.

As before, let
$f(u,n)=\prod_{i=1}^ng(u_i)$ where $0\leq g\leq1$ and $g(u_i)=1$ for
$u_i\geq u_g$.
Let
\[
A_rf(u,n)=f(u,n)\sum_{i=1}^n2a\int_{u_i}^{\infty}e^{-v/r}\bigl(g(v)-
1\bigr)\,dv+f(u,n)\sum_{i=1}^nG_r(u_i)\frac{g'(u_i)}{g(u_i)},
\]
where $G_r$ will be determined below. Note that a
particle at level $u_i$ is giving birth at rate $2are^{-u_i/r}$, and
the levels satisfy
\[
\dot{U}_i(t)=G_r(U_i(t)).
\]
Let $\alpha_r(n,du)$ be the distribution of $n$ independent
exponential random variables with mean $r$, and
define $e^{-\lambda_g}=r^{-1}\int_0^{\infty}g(v)e^{-v/r}\,dv$ so
\[
\widehat{f}(n)=\int f(u,n)\alpha_r(n,du)=e^{-\lambda_gn}.
\]
To calculate
$\int A_rf(u,n)\alpha_r(n,du)$, observe that
\begin{eqnarray*}
&&r^{-1}2a\int_0^{\infty}e^{-z/r}g(z)\int_z^{\infty}e^{-v/r}\bigl(g(v)-1
\bigr)\,dv\,dz\\
&&\qquad=are^{-2\lambda_g}-2a\int_0^{\infty}e^{-2z/r}g(z)\,dz,
\end{eqnarray*}
and assuming $G_r(0)=0$,
\begin{eqnarray*}
&&
r^{-1}\int_0^{\infty}e^{-z/r}G_r(z)g'(z)\,dz\\
&&\qquad=-r^{-1}\int_0^{\infty}
e^{-z/r}\bigl(G'_r(z)-r^{-1}G_r(z)\bigr)\bigl(g(z)-1\bigr)\,dz\\
&&\qquad=-r^{-1}\int_0^{\infty}e^{-z/r}\bigl(G'_r(z)-r^{-1}G_r(z)\bigr)g(z)\,dz\\
&&\qquad\quad{} +r^{-1}\int_0^{\infty}e^{-z/r}\bigl(G'_r(z)-r^{-1}G_r(z)\bigr)\,dz.
\end{eqnarray*}
Then for
\[
G'_r(z)-r^{-1}G_r(z)=e^{z/r}\frac d{dz}(e^{-z/r}G_r(z))=2ar(1-e^{
-z/r})-b,
\]
we have
\[
e^{-z/r}G_r(z)=2ar\biggl(r(1-e^{-z/r})-\frac r2(1-e^{-2z/r})\biggr)-br(1-e^{
-z/r})
\]
and
\begin{eqnarray*}
&&\int A_rf(u,n)\alpha_r(n,du)\\
&&\qquad =ne^{-\lambda_g(n-1)}\biggl(are^{-2\lambda_g}
-r^{-1}\int_0^{
\infty}e^{-z/r}\bigl(G'_r(z)-r^{-1}G_r(z)\\
&&\qquad\quad\hspace*{200pt}{}+2are^{-z/r}\bigr)g(z)\,dz\\
&&\qquad\quad\hspace*{113.4pt}{} +r^{-1}\int_0^{\infty}e^{-z/r}\bigl(G'_r(z)-r^{-1}
G_r(z)\bigr)\,dz\biggr)\\
&&\qquad =C_r\widehat{f}(n),
\end{eqnarray*}
where
%
%e4.10 ###
\begin{equation}\label{brproj2}\quad
C_r\widehat{f}(n)=ran\bigl(\widehat{f}(n+1)-\widehat{f}(n)\bigr)+(ra
-b)n\bigl(\widehat{f}(n-1)-\widehat{f}(n)\bigr)
\end{equation}
is the generator of a branching process.

Note that as $r\rightarrow\infty$, $G_r(z)$ converges to $az^2-bz$, and hence,
$A_r$ converges to $A$ given by (\ref{felldiff}).

%s4.6 ###
\subsection{Multitype branching processes} We now consider
a branching particle system with $m$ possible types,
$S=\{1,2,\ldots,m\}$.
We assume that a particle of type $\zeta_1\in S$ gives birth
to a particle of type $\zeta_2\in S$ at rate $ra^{(r)}(\zeta_1,\zeta_
2)$ and dies
at rate $ra^{(r)}(\zeta_1)-b^{(r)}(\zeta_1)$, where $a^{(r)}(\zeta_
1)=\sum_{j\in S}a^{(r)}(\zeta_1,j)$.

The fact that the ordered representations constructed
for the previous examples give the correct
measure-valued processes depends on the fact that
observing a birth event in the measure-valued process
gives no information about the levels of the particles
after the birth event. That, in turn, depends on the
offspring being indistinguishable from the parent. Since
in the current model, the type of an offspring may
differ from the type of the parent, we need to find a
way to ``preserve ignorance'' about the levels when the
type of the offspring is different. We accomplish this
goal by
randomizing the assignment of the parent and offspring to
the original level of the parent and a new level. Let
$f(\zeta,u,n)$ be of the form
\[
f(\zeta,u,n)=\prod_{i=1}^ng(\zeta_i,u_i).
\]
Then the generator of the ordered representation of the
branching process described above is given by
\begin{eqnarray*}
&&A_rf(\zeta,u,n)\\
&&\qquad =f(\zeta,u,n) \sum_{i=1}^n \sum_{j\in S} 2 a^{(r)}(\zeta_
i,j)\\
&&\qquad\quad\hspace*{74pt}{}\times\int_{u_i}^r\biggl[\frac12\biggl(\frac{g(\zeta_i,u_i) g(j,v)
+g(\zeta_i,v) g(j,u_i)}{g(\zeta_i,u_i)}\biggr)-1\biggr] \,dv\\
&&\qquad\quad{} +f(\zeta,u,n) \sum_{i=1}^n \bigl[a^{(r)}(\zeta_
i) u_i^2-b^{(r)}(\zeta_i) u_i\bigr]\,\frac{\partial_{u_i}g(\zeta_
i,u_i)}{g(\zeta_i,u_i)},
\end{eqnarray*}
where as before, each level satisfies
\[
\frac d{dt}U^{(r)}_i(t)=a^{(r)}(X_i(t))U_i^2(t)-b^{(r)}(X_i(t))
U_i(t).
\]
Let $Q^{(r)}h(\zeta)=\sum_{j\in S}a^{(r)}(\zeta,j)[h(j)-h(
\zeta)]$.
Because of the randomization of the level assignments at
each birth event, it follows that
\[
h(X_i(t))-\int_{\tau_i}^t\bigl(r-U_i(s)\bigr)Q^{(r)}h(X_i(s))\,ds
\]
is a martingale.

Taking $\alpha_{r}(n,du)$ as before, we have that
\begin{eqnarray*}
\int A_rf(\zeta,u,n) \alpha(n,du)&=&C_r\overline{f}(\zeta,n)\\
&=&\overline f(\zeta,n)\sum_{i=1}^n\sum_{j\in S}r a^{(r)}(\zeta_
i,j)[\overline g(j)-1]\\
&&{} +\overline f(\zeta,n)\sum_{i=1}^n\bigl[r a^{(r)}(\zeta_
i)-b^{(r)}(\zeta_i)\bigr]\biggl[\frac1{\overline{g}(\zeta_i)}-1
\biggr],
\end{eqnarray*}
where $\overline f(\zeta,n)=\prod_{i=1}^n\overline g(\zeta_i)$ and
$\overline
g(\zeta_i)=\frac1r\int_0^rg(\zeta_i,v) \,dv$.
Hence, $C_r\overline{f}(\zeta,n)$ is the generator of a multitype
branching process.

Assume that
\begin{eqnarray*}
a(\zeta,j)&=&\lim_{r\rightarrow\infty}a^{(r)}(\zeta,j),\\
a(\zeta)&=&\lim_{r\rightarrow\infty}a^{(r)}(\zeta)=\sum_{j\in S}
a(\zeta,j),\\
b(\zeta)&=&\lim_{r\rightarrow\infty}b^{(r)}(\zeta)
\end{eqnarray*}
and that
\[
Qh(\zeta)=\sum_{j\in S}a(\zeta,j)[h(j)-h(\zeta)]
\]
is the generator of an irreducible, finite
state Markov chain. Let $\pi$ denote the unique
stationary distribution for $Q$. It is clear from the ergodicity of
the Markov chain that in the limit, the levels must
satisfy
\[
\frac d{dt}U_i(t)=\overline{a}U_i^2(t)-\overline{b}U_i(t),
\]
where
$\overline{a}=\sum_j\pi(j)a(j)$ and $\overline{b}=\sum_j\pi(j)b(j)$.

We can make this observation precise by analyzing the
asymptotic behavior of the generator.
Taking $g(\zeta,u)=\exp(-h_0(u)+\frac1rh(\zeta,u))$, where
$h(\zeta,u)$ and $h_0(u)$ are equal to zero if $u\geq u_g$, and
letting $r\rightarrow\infty$, we have that
\begin{eqnarray*}
\lim_{r\rightarrow\infty}f(\zeta,u)&=&\lim_{r\rightarrow\infty}\exp
\biggl(-\sum_ih_0(u_i)+\frac1r\sum_ih(\zeta_i,u_i)\biggr)\\
&=&\exp\biggl(-\sum_
ih_0(u_i)\biggr)\equiv\overline{f}(u)
\end{eqnarray*}
and since
\begin{eqnarray}
&&\frac{g(\zeta_i,u_i) g(j,v)+g(\zeta_i,v) g(j,u_i)}{g(\zeta_i
,u_i)}\nonumber\\
&&\qquad=e^{-h_0(v)}\bigl(e^{r^{-1}h(j,v)}+e^{r^{-1}(h(\zeta_i,v)+h
(j,u_i)-h(\zeta_i,u_i))}\bigr),
\nonumber\\
%
%e4.11 ###
\label{multlim}
&&\lim_{r\rightarrow\infty}A_rf(\zeta,u)\nonumber\\
&&\qquad=\overline f(u)\sum_{i }
\biggl\{2 a(\zeta_i)\int_{u_i}^{u_g}\bigl[e^{-h_0(v)}-1\bigr] \,dv
\nonumber\\[-8pt]\\[-8pt]
&&\qquad\quad\hspace*{40.8pt}{} +\sum_{j\in S}a(\zeta_i,j)[h(j,u_i)-h(\zeta_
i,u_i)]\nonumber\\
&&\qquad\quad\hspace*{54.1pt}{} -[a(\zeta_i)u_i^2-b(\zeta_i)u_i
] \,\partial_{u_i}h_0(u_i)\biggr\}.\nonumber
\end{eqnarray}

If $\sum\pi(j) G(j,u)\equiv0$ for all $u$, then there
exists $h$ such that
\[
\sum_{j\in S}a(\zeta,j)[h(j,u)-h(\zeta,u)]= G(\zeta
,u).
\]
Consequently, there exists $h$ such that the right-hand side of
(\ref{multlim}) becomes
\[
A\overline{f}(u)=\overline f(u)\sum_{i }\biggl\{2 \overline{a}\int_{u_i}^{
u_g}\bigl[e^{-h_0(v)}-1\bigr] \,dv-[\overline au_i^2-\overline bu_i
] \,\partial_{u_i}h_0(u_i)\biggr\},
\]
which is just a rewriting of (\ref{felldiff}), and hence
we have convergence of the normalized total population
to the Feller diffusion.

For earlier work, see \cite{GL90,JM86,Kur78} and Section 9.2
of \cite{EK86}.

%s4.7 ###
\subsection{Models with catastrophic death}
Now consider
%
%e4.12 ###
\begin{eqnarray}\label{fingencat}
A_rf(x,u,n)&=&f(x,u,n)\sum_{i=1}^n\frac{Bg(x_i,u_i)}{g(x_i,u_i)}\nonumber\\
&&{}
+f(x,u,n)\sum_{i=1}^n2a(x_i)\int_{u_i}^r\bigl(g(x_i,v)-1\bigr)\,dv\nonumber\\[-8pt]\\[-8pt]
&&{} +f(x,u,n)\sum_{i=1}^n\bigl(a(x_i)u_i^2-b(x_i)u_i\bigr)\,\frac{\partial_{
u_i}g(x_i,u_i)}{g(x_i,u_i)}\nonumber\\
&&{} +\int_V\bigl(f(x,c(u,x,v),n)-f(x,u,n)\bigr)\gamma(dv),\nonumber
\end{eqnarray}
where $\gamma$ is a $\sigma$-finite measure on a measurable space
$(V,\mathcal{V})$,
\[
c(u,x,v)=(u_1\rho(x_1,v),u_2\rho(x_2,v),\ldots)
\]
and $\rho(x_i,v)\geq1$. Then as in Section \ref{sec-bran}
\begin{eqnarray*}
C_r\widehat{f}(x,n)&=&\sum_{i=1}^nB_{x_i}\widehat{f}(x,n)+\sum_{i=1}^nr
a(x_i)\bigl(\widehat{f}\bigl(b(x|x_i),n+1\bigr)-\widehat{f}(x,n)\bigr)\\
&&{} +\sum_{i=1}^n\bigl(ra(x_i)-b(x_i)\bigr)\bigl(\widehat{f}\bigl(d(x|x_i),n-1\bigr)-\widehat{
f}(x,n)\bigr)\\
&&{} +\int_V\Biggl(\prod_{i=1}^n\bigl(\rho^{-1}(x_i,v)\widehat g(x_i)+\bigl(
1-\rho^{-1}(x_i,v)\bigr)\bigr)-\widehat{f}(x,n)\Biggr)\gamma(dv).
\end{eqnarray*}
For simplicity, assume that $\gamma(V)<\infty$. Then at rate $\gamma
(V)$
an event occurs in which an element $v$ is selected from
$V$, and given $v$, particles are independently killed, with
the probability that a particle at $x_i$ survives being $\rho^{-1}
(x_i,v)$.

Letting $r\rightarrow\infty$ to obtain $A$ and integrating,
\begin{eqnarray*}
\alpha Af&=&\int_E\bigl(-Bh(y)+a(y)h(y)^2-b(y)h(y)\bigr)\mu(dy)\exp
\{-\langle h,\mu\rangle\}\\
&&{} +\int_V\bigl(\exp\{-\langle\rho^{-1}(\cdot,v)h,\mu\rangle\}-\exp
\{-\langle h,\mu\rangle\}\bigr)\gamma(dv)\\
&=&C\widehat{f}(\mu).
\end{eqnarray*}

Branching processes with catastrophes have been
considered in a series of papers by Pakes
\cite{Pak86,Pak87,Pak88,Pak89a,Pak89b,Pak90} and by Grey
\cite{Gre88}.

\begin{appendix}\label{app}
\section*{Appendix}

%s4.8 ###
\subsection{Poisson random measures} Let $(S,\mathcal{S})$ be a
measurable space, and let $\nu$ be a $\sigma$-finite measure on $
\mathcal{S}$.
$\xi$ is a \textit{Poisson random measure} with mean measure $\nu$ if:

\begin{enumerate}[(a)]
\item[(a)]$\xi$ is a random counting
measure on $S$;

\item[(b)]for each $A\in\mathcal{S}$ with $\nu(A)<\infty$, $\xi(A)$ is
Poisson distributed with parameter $\nu(A)$;

\item[(c)]for
$A_1,A_2,\ldots\in\mathcal{S}$ disjoint, $\xi(A_1),\xi(A_2),\ldots$
are independent.
\end{enumerate}
\setcounter{theorem}{0}
\begin{lemma}\label{poistran}
If $H\dvtx S\rightarrow S_0$ is Borel measurable and $\widehat{\xi}(A)=\xi
(H^{-1}(A))$, then $\widehat{\xi}$ is a Poisson
random measure on $S_0$ with mean measure $\widehat{\nu}$ given by
$\widehat{\nu}(A)=\nu(H^{-1}(A))$.
\end{lemma}
\begin{remark}
$\widehat{\nu}$ need not be $\sigma$-finite even if $\nu$ is, but the
meaning of
the lemma should still be clear. $\sigma$-finite or not,
$\widehat{\nu}(A)=\infty$ if and only if $\widehat{\xi}(A)=\infty$ a.s.
\end{remark}
\begin{pf*}{Proof of Lemma \ref{poistran}}
The lemma follows from the fact
that $A_1,A_2,\ldots$ disjoint implies $H^{-1}(A_1),H^{-1}(A_2),\ldots$ are
disjoint.
\end{pf*}
\begin{lemma}\label{lemma2}
If $\xi$ is a Poisson random measure with mean measure $\nu$
and $f\in L^1(\nu)$, then
%
%e4.13 ###
\setcounter{equation}{0}
\begin{eqnarray}\label{lapfunc}
E\bigl[e^{\int f(z)\xi(dz)}\bigr]&=&e^{\int(e^f-1)\,d\nu},
\\
%
%e4.14 ###
\label{poisvar}
E\biggl[\int f(z)\xi(dz)\biggr]&=&\int f\,d\nu,\qquad
\operatorname{Var}\biggl(\int f(z)\xi(dz)\biggr)=\int f^2\,d\nu,
\end{eqnarray}
allowing $\infty=\infty$.

Letting $\xi=\sum_i\delta_{Z_i}$, for $g\geq0$ with
$\log g\in L^1(\nu)$,
\[
E\biggl[\prod_ig(Z_i)\biggr]=e^{\int(g-1)\,d\nu}.
\]
Similarly, if $hg,g-1\in L^1(\nu)$, then
\[
E\biggl[\sum_jh(Z_j)\prod_ig(Z_i)\biggr]=\int hg\,d\nu\, e^{\int(g-1)\,d\nu}
\]
and
\[
E\biggl[\sum_{i\neq j}h(Z_i)h(Z_j)\prod_kg(Z_k)\biggr]=\biggl(\int hg\,d\nu\biggr)^2e^{\int
(g-1)\,d\nu}.
\]
\end{lemma}
\begin{pf}
The independence properties of $\xi$ imply (\ref{lapfunc}) for
simple functions. The general case follows by
approximation. The other identities follow in a similar
manner. Note that the integrability of the random
variables in the expectations above can be verified by
replacing $g$ by $(|g|\wedge a){\mathbf1}_A+{\mathbf1}_{A^c}$ and $h$ by $
(|h|\wedge a){\mathbf1}_A$ for
$0<a<\infty$ and $\nu(A)<\infty$ and passing to the limit as
$a\rightarrow\infty$ and $A\nearrow E$.
\end{pf}
\begin{lemma}\label{poisscale}
If $\xi_0=\sum_i\delta_{U_i}$ is a Poisson random measure on $[0
,\infty)$ with mean
measure $\lambda\varLambda$, $\varLambda$ Lebesgue measure, and $
\{X_i\}$ are i.i.d.
positive random variables, independent of $\xi_0$, then
\[
\xi=\sum_i\delta_{(X_i,U_i)}
\]
is a Poisson random measure on $[0,\infty)^2$
with mean measure $\lambda\mu_X\times\varLambda$, were $\mu_{X}$ is the law
of $X_{1}$.

If $\kappa=E[\frac1{X_i}]<\infty$, then
\[
\widehat{\xi}=\sum_{}\delta_{X_iU_i}
\]
is a Poisson random measure on $[0,\infty)$ with
mean measure $\lambda\kappa\varLambda$.
\end{lemma}
\begin{pf}
By Lemma \ref{poistran}, $\widehat{\xi}$ is a Poisson random measure
with mean measure given by
\begin{eqnarray*}
\widehat{\nu}[0,c]&=&\lambda\mu_X\times\varLambda\{(x,u)\dvtx xu\leq c\}
=\lambda\int_0^{\infty}P\{X^{-1}\geq uc^{-1}\}\,du\\
&=&\lambda cE[X^{-1}].
\end{eqnarray*}
\upqed
\end{pf}

%s4.9 ###
\subsection{Conditionally Poisson systems}\label{sectcp}

We begin by considering general conditionally Poisson
systems or Cox processes. Consider $(S,d)$ a metric
space, and let $\xi$ be a random counting measure on $S$ and
$\Xi$ be a locally finite random measure on $S$. [A measure
$\nu$ on $S$ is locally finite if for each $x\in S$, there exists an
$\epsilon>0$ such that $\nu(B_{\varepsilon}(x))<\infty$.] We say that $
\xi$ is
conditionally Poisson with Cox measure $\Xi$ if, conditioned
on $\Xi$, $\xi$ is a Poisson random measure with mean
measure~$\Xi$. This requirement is equivalent to
\[
E[e^{-\int_Sf\,d\xi}]=E\bigl[e^{-\int_S(1-e^{-f})\,d\Xi}\bigr],
\]
for all nonnegative $f\in M(S)$, where $M(S)$ is the set of all Borel
measurable functions on $S$. Since the collection of
functions $F_f(\mu)=e^{-\int_Sf\,d\mu}$ is closed under multiplication
and separates points in the space of locally finite
measures, the distribution of $\Xi$ determines the
distribution of $\xi$.

We are actually interested in the conditionally Poisson
system on $S\times[0,\infty)$ with Cox measure $\Xi\times\varLambda$,
where $
\varLambda$ is
Lebesgue measure. Then for nonnegative $f\in M(S)$, we
have
\[
E[e^{-\int_{S\times[0,K]}f\,d\xi}]=E\bigl[e^{-K\int_S(1-e^{-f})\,d\Xi}\bigr],
\]
and the distribution of $\xi$ determines the distribution of
$\Xi$, where we consider $f\in M(S)$ to be a function on
$S\times[0,K]$ satisfying $f(x,u)=f(x)$.
In particular,
\[
\Xi(f)=\lim_{K\rightarrow\infty}\frac1K\int_{S\times[0,K]}f\,d
\xi  \qquad\mbox{a.s.}
\]
\begin{lemma}\label{stbnd}
Suppose $\xi$ is a conditionally Poisson random measure on
$S\times[0,\infty)$ with Cox measure $\Xi\times\varLambda$, and let $
f\in M(S)$,
$0\leq f\leq1$. Then for $C,D>0$,
%
%e4.15 ###
\begin{equation}\label{xibnd}
P\biggl\{\int_{S\times[0,K]}f\,d\xi\geq C\biggr\}\leq\frac{KD}
C+P\biggl\{\int_Sf\,d\Xi\geq D\biggr\}
\end{equation}
and
%
%e4.16 ###
\begin{equation}\label{Xibnd}
P\biggl\{\int_Sf\,d\Xi\geq C\biggr\}\leq\frac{E[1-e^{-C^{-1}\int_{
S\times[0,K]}f\,d\xi}]}{1-e^{-Ke^{-C^{-1}}}}.
\end{equation}

Let $\{\xi_{\alpha},\alpha\in\mathcal{A}\}$ be a collection of
conditionally Poisson
random measures on $S\times[0,\infty)$ with Cox measures
$\Xi_{\alpha}\times\varLambda$,
and let $f\in M(S)$,
$0\leq f\leq1$. Then $\{\int_{S\times[0,K]}f\,d\xi_{\alpha},\alpha
\in\mathcal{A}\}$ is stochastically
bounded if and only if $\{\int_Sf\,d\Xi_{\alpha},\alpha\in\mathcal{A}
\}$ is stochastically
bounded.
\end{lemma}
\begin{pf}
Since $E[\int_{S\times[0,K]}f\,d\xi|\Xi]=K\int_Sf\,d\Xi$,
\[
P\biggl\{\int_{S\times[0,K]}f\,d\xi\geq C\Big|\Xi\biggr\}\leq\frac{K\int_Sf\,d\Xi}
C\wedge1\leq\frac{KD}C+{\mathbf1}_{\{\int_Sf\,d\Xi\geq D\}},
\]
and taking expectations gives (\ref{xibnd}).

By (\ref{lapfunc})
\begin{eqnarray*}
E[1-e^{-\int_{S\times[0,K]}\varepsilon f\,d\xi}]&=&E\bigl[1-e^{-K\int_S(1-
e^{-\varepsilon f})\,d\Xi}\bigr]\\
&\geq&E[1-e^{-\varepsilon Ke^{-\varepsilon}\int_Sf\,d\Xi}]\\
&\geq&(1-e^{-\varepsilon Ke^{-\varepsilon}C})P\biggl\{\int_Sf\,d\Xi\geq C\biggr\},
\end{eqnarray*}
and taking $\varepsilon=C^{-1}$ gives (\ref{Xibnd}).

The final statement follows from the two inequalities.
\end{pf}

Let $\widehat{\xi}=\sum\delta_{X_i}$ be a point process on $S$, and let $
\{U_i\}$ be
independent random variables, uniformly distributed on
$[0,r]$ and independent of $\widehat{\xi}$. Define
%
%e4.17 ###
\begin{equation}\label{ulev}
\xi=\sum\delta_{(X_i,U_i)},\qquad \Xi_r=r^{-1}\widehat{
\xi}.
\end{equation}
Then for $f\geq0$ on $S\times[0,r]$,
%
%e4.18 ###
\begin{equation}\label{uniid}
E[e^{-\int_{S\times[0,r]}f\,d\xi}|\Xi_r]=\prod_i
\biggl(r^{-1}\int_0^re^{-f(X_i,u)}\,du\biggr)=e^{-\int_SF_f^r(x)\Xi_r(dx)},\hspace*{-28pt}
\end{equation}
where
\[
F_f^r(x)=-r\log\frac1r\int_0^re^{-f(x,u)}\,du=-r\log\biggl(1-\frac1r\int_
0^r\bigl(1-e^{-f(x,u)}\bigr)\,du\biggr).
\]
We have the following analog of Lemma \ref{stbnd}.
\begin{lemma}\label{stbndun}
Suppose $\xi$ and $\Xi_r$ are given by (\ref{ulev}),
and let $f\in M(S)$,
$0\leq f\leq1$. Then for $C,D>0$ and $K \leq r$,
%
%e4.19 ###
\begin{equation}\label{xibndun}
P\biggl\{\int_{S\times[0,K]}f\,d\xi\geq C\biggr\}\leq\frac{KD}
C+P\biggl\{\int_Sf\,d\Xi_r\geq D\biggr\}
\end{equation}
and
%
%e4.20 ###
\begin{equation}\label{Xibndun}
P\biggl\{\int_Sf\,d\Xi_r\geq C\biggr\}\leq\frac{E[1-e^{-C^{-1}
\int_{S\times[0,K]}f\,d\xi}]}{1-e^{-Ke^{-C^{-1}}}}.
\end{equation}
\end{lemma}
\begin{pf}
Since $E[\int_{S\times[0,K]}f\,d\xi|\Xi_r]=K\int_Sf\,d\Xi_r$,
\[
P\biggl\{\int_{S\times[0,K]}f\,d\xi\geq C\Big|\Xi_r\biggr\}\leq\frac{K\int_Sf\,d\Xi_
r}C\wedge1\leq\frac{KD}C+{\mathbf1}_{\{\int_Sf\,d\Xi_r\geq D\}},
\]
and taking expectations gives (\ref{xibndun}).

Defining
\[
G_{K,\varepsilon,f}^r(x)=-r\log\biggl(1-\frac Kr\bigl(1-e^{-\varepsilon f(x)}\bigr)\biggr)\geq
\varepsilon Ke^{-\varepsilon}f(x),
\]
by (\ref{uniid})
\begin{eqnarray*}
E[1-e^{-\int_{S\times[0,K]}\varepsilon f\,d\xi}]&=&E[1-e^{-\int_SG^r_{
K,\varepsilon,f}\,d\Xi_r}]\\
&\geq&E[1-e^{-\varepsilon Ke^{-\varepsilon}\int_Sf\,d\Xi_r}]\\
&\geq&(1-e^{-\varepsilon Ke^{-\varepsilon}C})P\biggl\{\int_Sf\,d\Xi_r\geq C\biggr\},
\end{eqnarray*}
and taking $\varepsilon=C^{-1}$ gives (\ref{Xibndun}).
\end{pf}
\begin{lemma}
Suppose $\xi$ is a conditionally Poisson random measure on
$S\times[0,\infty)$ with Cox measure $\Xi\times\varLambda$.
If $\Xi(S) <\infty$ a.s., then we can write $\xi=\sum_{i=1}^{\infty}
\delta_{(X_i,U_i)}$ with
$U_1<U_2<\cdots$ a.s. and $\{X_i\}$ exchangeable.
\end{lemma}
\begin{pf}
Let $\{\widetilde{X}_i\}$ be exchangeable\vspace*{1pt} with de Finetti measure $\frac{
\Xi}{|\Xi|}$, and
let $Y$ be a unit Poisson process with jump times $\{S_i\}$
independent of of $\{\widetilde{X}_i\}$ and $\Xi$.
Define $\widetilde{\xi}=\sum_{i=1}^{\infty}\delta_{(\widetilde{X}_i,|\Xi
|^{-1}S_i)}$,
and note that
\begin{eqnarray*}
E[e^{-\int fd\widetilde{\xi}}]&=&E\biggl[\prod_ie^{-f(\widetilde{X}_i,|\Xi|^{
-1}S_i)}\biggr]\\
&=&E\biggl[\prod_i\int e^{-f(z,|\Xi|^{-1}S_i)}|\Xi|^{-1}\Xi(dz)\biggr]\\
&=&E\biggl[\exp\biggl\{-\int_0^{\infty}\biggl(1-\int e^{-f(z,|\Xi|^{-1}s)}|\Xi|^{
-1}\Xi(dz)\biggr)\,ds\biggr\}\biggr]\\
&=&E\biggl[\exp\biggl\{-\int_0^{\infty}\int\bigl(1- e^{-f(z,u)}\bigr) \Xi(dz)\,du\biggr\}\biggr].
\end{eqnarray*}
Consequently, $\widetilde{\xi}$ is conditionally Poisson with Cox
measure $\Xi\times\varLambda$, and $\xi$ and $\widetilde{\xi}$ have the
same distribution.
\end{pf}

As in Lemma \ref{poisscale}, we have the following.
\begin{lemma}
Suppose $\xi$ is a conditionally Poisson random measure on
$S\times[0,\infty)^2$ with Cox measure $\Xi\times\varLambda$, where $
\Xi$
is a random measure on $S\times[0,\infty)$.
Suppose
\[
\widehat{\Xi}(A)=\int_{S\times[0,\infty)}\frac1y{\mathbf1}_A(x)\Xi
(dx\times dy)
\]
defines a locally finite random measure on $S$.
Then writing
$\xi=\sum_i\delta_{(X_i,Y_i,U_i)}$, $\widehat{\xi}=\sum_i\delta_{(X_
i,Y_iU_i)}$ is a conditionally\vspace*{1pt}
Poisson random measure on $S\times[0,\infty)$ with Cox measure
$\widehat{\Xi}\times\varLambda$,
and hence
\[
\int_{S\times[0,\infty)}y^{-1}f(x)\Xi(dx\times dy)=\widehat{\Xi}
(f)=\lim_{K\rightarrow\infty}\frac1K\int_{S\times[0,K]}f\,d\widehat{\xi}
\qquad\mbox{a.s.}
\]
\end{lemma}

%s4.10 ###
\subsection{Convergence results}\label{convres}
Let $\{h_k,k=1,2,\ldots\}\subset\overline{C}(S)$ satisfy $0\leq h_k\leq
1$ and
$\bigcup_k\{x\dvtx h_k(x)>0\}=S$, where $\overline{C}(S)$ denotes the space of
bounded continuous functions on $S$, and let
$\mathcal{M}_{\{h_k\}}(S)$ be the collection of Borel measures on $S$
satisfying $\int_Sh_k\,d\nu<\infty$, for all $k$, topologized by the requirement
that $\nu_n\rightarrow\nu$ if and only if $\int_Sfh_k\,d\nu_n\rightarrow
\int_Sfh_k\,d\nu$ for all
$f\in\overline{C}(S)$ and $k$; that is, the measures $d\nu_n^k=h_k\,d\nu_
n$ converge
weakly for each $k$.
Similarly, let $\mathcal{M}_{\{h_k\}}(S\times[0,\infty))$ be the
space of Borel measures on $S\times[0,\infty)$ satisfying
$\int_{S\times[0,K]}h_k\,d\mu<\infty$ for all $k=1,2,\ldots$ and $
K>0$, topologized by
the requirement that $\mu_n\rightarrow\mu$ if and only if
\[
\int_{S\times[0,\infty)}fh_k\,d\mu_n\rightarrow\int_{S\times[0
,\infty)}fh_k\,d\mu,
\]
for all $k$ and $f\in\overline{C}(S\times[0,\infty))$ such that the support
of $f$ is contained in $S\times[0,K]$ for some $K>0$. Note that
in both cases,
$\mathcal{M}_{\{h_k\}}$ is metrizable. To simplify notation, let
\[
\mathcal{C}_{\{h_k\}}(S)=\{f\in\overline{C}(S)\dvtx|f|\leq ch_k\mbox{ for
some }
c>0\mbox{ and }h_k\}.
\]
Then convergence in $\mathcal{M}_{\{h_k\}}(S)$ is equivalent to
convergence of $\int_Sf\,d\nu_n$ for all $f\in\mathcal{C}_{\{h_k\}}(S
)$.
\begin{theorem}\label{convequiv1}
Let $\{\xi^n\}$ be a sequence of conditionally Poisson random
measures on $S\times[0,\infty)$ with Cox measures $\{\Xi^n\times
\varLambda\}$. Then
$\xi^n\Rightarrow\xi$ in $\mathcal{M}_{\{h_k\}}(S\times[0,\infty))$
if and only if $
\Xi^n\Rightarrow\Xi$ in
$\mathcal{M}_{\{h_k\}}(S)$. If the limit holds, then $\xi$ is conditionally
Poisson with Cox measure $\Xi\times\varLambda$.
\end{theorem}
\begin{pf}
Suppose $\xi^n\Rightarrow\xi$ in $\mathcal{M}_{\{h_k\}}(S\times
[0,\infty
))$. Then for each $f\in\overline{C}(S)$,
$f\geq0$,
each $k$, and all but
countably many $K$
\begin{eqnarray*}
E[e^{-\int_{S\times[0,K]}fh_k\,d\xi}]&=&\lim_{n\rightarrow\infty}E
[e^{-\int_{S\times[0,K]}fh_k\,d\xi^n}]\\
&=&\lim_{n\rightarrow\infty}E\bigl[
e^{-K\int_Sh_k^{-1}(1-e^{-fh_k})h_k\,d\Xi^n}\bigr].
\end{eqnarray*}
For $g\geq0$ and $K$ satisfying $\sup_xK^{-1}g(x)h_k(x)<1$, let
\[
f(x)=\cases{
-h_k^{-1}(x)\log\bigl(1-K^{-1}g(x)h_k(x)\bigr), &\quad  $h_k(x)>0$,\vspace*{1pt}\cr
K^{-1}g(x), &\quad  $h_k(x)=0$,}
\]
and we see that
\[
\lim_{n\rightarrow\infty}E[e^{-\int_Sgh_k\,d\Xi^n}]=E[e^{-\int_{S
\times[0,K]}fh_k\,d\xi}]
\]
exists. Since $\xi^n\Rightarrow\xi$ in $\mathcal{M}_{\{h_k\}}(S\times
[0,\infty))$, $\{\int_{S\times[0,K]}h_k\,d\xi^n\}$ is
stochastically bounded and by Lemma \ref{stbnd},
$\{\int h_k\,d\Xi^n\}$ must be stochastically bounded.
Tightness follows similarly. Consequently, $\{\Xi^n\}$ is
relatively compact in $\mathcal{M}_{\{h_k\}}(S)$ in distribution, and the
unique limit $\Xi$ is determined by the fact that
\[
E[e^{-\int_Sgh_k\,d\Xi}]=E[e^{-\int_{S\times[0,K]}fh_k\,d\xi}],
\]
for $g$ and $f$ related as above. The proof of the converse
is similar.
\end{pf}
\begin{theorem}\label{convequiv2}
For each $n=1,2,\ldots,$ let $r_n>0$ and $\xi^n$ be a point process on
$S\times[0,r_n]$, and define
%
%e4.21 ###
\begin{equation}\label{finmeas}
\Xi^n(dx)=\frac1{r_n}\xi^n(dx\times[0,r_n]).
\end{equation}
Suppose for $f\geq0$,
$E[e^{-\int f(x,u)\xi^n(dx\times du)}]=E[e^{-\int F^n_f(x)\Xi^n(d
x)}]$,
where
\[
F^n_f(x)=-r_n\log\frac1{r_n}\int_0^{r_n}e^{-f(x,u)}\,du=-r_n\log
\biggl(1-\frac1{r_n}\int_0^{r_n}\bigl(1-e^{-f(x,u)}\bigr)\,du\biggr),
\]
that is, the $[0,r_n]$ components are independent, uniformly
distributed, and independent of $\Xi^n$.
Then assuming $r_n\rightarrow\infty$,
$\xi^n\Rightarrow\xi$ in $\mathcal{M}_{\{h_k\}}(S\times[0,\infty))$
if and only if $
\Xi^n\Rightarrow\Xi$ in
$\mathcal{M}_{\{h_k\}}(S)$.
If the limit holds, then $\xi$ is conditionally
Poisson with Cox measure $\Xi\times\varLambda$.
\end{theorem}
\begin{pf}
For $g_0,f_0\geq0$, $g_0\in C_c([0,\infty))$, $f_0\in\overline{C}(
S)$
and
$f(x,u)=h_k(x)\times f_0(x)g_0(u)$, $F^n_f(x)\rightarrow\int_0^{\infty}(
1-e^{-f(x,u)})\,du$,
and\vspace*{1pt} the remainder of the proof is similar to that of Theorem
\ref{convequiv1}.
\end{pf}

These convergence theorems apply only to the one-dimensional
distributions of the models considered in this
paper.
To address convergence as processes, note that for finite
$r$ and $\Xi_r(t,dx)=r^{-1}\xi(t,dx\times[0,r])$, the models satisfy
%
%e4.22 ###
\begin{equation}\label{cndrelfi}
E\bigl[e^{-\int_{S\times[0,r]}f(x,u)\xi(t,dx\times d
u)}|\mathcal{F}_t^{\Xi_r}\bigr]=e^{-\int F^r_f(x)\Xi_r(t,dx)},
\end{equation}
where
\[
F_f^r(x)=-r\log\frac1r\int_0^re^{-f(x,u)}\,du=-r\log\biggl(1-\frac1r\int_
0^r\bigl(1-e^{-f(x,u)}\bigr)\,du\biggr),
\]
and the $r=\infty$ models satisfy
%
%e4.23 ###
\begin{equation}\label{cndrel}
E\bigl[e^{-\int_{S\times[0,K]}f\,d\xi(t)}|\mathcal{F}_t^{
\Xi}\bigr]=e^{-K\int_S(1-e^{-f})\,d\Xi(t)},
\end{equation}
for $f\in\mathcal{C}_{\{h_k\}}(S)$. The following estimates imply that
convergence of the finite-dimensional distributions for $\xi^n$
imply convergence of the finite-dimen\-sional distributions
for $\Xi^n$ (or $\Xi^n_{r_n}$, assuming $r_n\rightarrow\infty$);
however, convergence
of the finite-dimensional distributions of $\Xi^n$ may not
imply convergence of the finite-dimensional distributions
of of $\xi^n$.
\begin{lemma}\label{coxappr}
Suppose $\xi$ is a conditionally Poisson random measure
on $S\times[0,\infty)$
with Cox\vadjust{\goodbreak} measure $\Xi\times\varLambda$, $\Xi$ with values in $\mathcal{M}_{
\{h_k\}}(S)$.
Then for each $f\in\mathcal{C}_{\{h_k\}}(S)$ and $\delta,K,K'>0$,
%
%e4.24 ###
\begin{eqnarray}\label{poisest}
&&
P\biggl\{\biggl|K^{-1}\int_{S\times[0,K]}f\,d\xi-\int_Sf\,d\Xi\biggr|\geq
\delta\biggr\}\nonumber\\
&&\qquad\leq\frac C{K\delta^2}+P\biggl\{\int f^2\,d\Xi>C\biggr\}
\\
&&\qquad\leq\frac C{K\delta^2}+\frac{E[1-e^{-C^{-1}\int_{S\times[0,K'
]}f^2\,d\xi}]}{1-e^{-K'e^{-C^{-1}}}}.\nonumber
\end{eqnarray}

Suppose $\xi$ satisfies (\ref{cndrelfi}).
Then for each $f\in\mathcal{C}_{\{h_k\}}(S)$, $\delta>0$ and $0<K
$, \mbox{$K'<r$}
%
%e4.25 ###
\begin{eqnarray}\label{uniest}
&&
P\biggl\{\biggl|K^{-1}\int_{S\times[0,K]}f\,d\xi-\int_Sf\,d\Xi_r\biggr|\geq
\delta\biggr\}\nonumber\\
&&\qquad\leq\frac{(r-K)C}{rK\delta^2}+P\biggl\{\int f^2\,d\Xi_r>C\biggr\}
\\
&&\qquad\leq\frac{(r-K)C}{rK\delta^2}+\frac{E[1-e^{-C^{-1}\int_{S\times
[0,K']}f^2\,d\xi}]}{1-e^{-K'e^{-C^{-1}}}}.\nonumber
\end{eqnarray}
\end{lemma}
\begin{pf}
By (\ref{poisvar}) and the Chebyshev inequality,
\[
P\biggl\{\biggl|K^{-1}\int_{S\times[0,K]}f\,d\xi-\int_Sf\,d\Xi\biggr|\geq
\delta\Big|\Xi\biggr\}\leq\frac{\int f^2\,d\Xi}{K\delta^2}\wedge1\leq\frac
C{K\delta^2}+{\mathbf1}_{\{\int f^2\,d\Xi>C\}},
\]
and taking expectations gives the first inequality in
(\ref{poisest}). The second inequality follows by
(\ref{Xibnd}).

Similarly, for the second part,
\begin{eqnarray*}
P\biggl\{\biggl|K^{-1}\int_{S\times[0,K]}f\,d\xi-\int_Sf\,d\Xi_r\biggr|
\geq\delta\Big|\Xi_r\biggr\}&\leq&\frac{\int(1-K/r)f^2\,d\Xi_r}{
K\delta^2}\wedge1\\
&\leq&\frac{(r-K)C}{rK\delta^2}+{\mathbf1}_{\{\int
f^2\,d\Xi_r>C\}},
\end{eqnarray*}
and taking expectations gives the first inequality in
(\ref{uniest}). The second inequality follows by
(\ref{Xibndun}).
\end{pf}

The estimates in Lemma \ref{coxappr} allow verifying
convergence of measure-valued processes satisfying
(\ref{cndrel}) or (\ref{cndrelfi}) by verifying convergence
of the corresponding particle representations.
\begin{theorem}\label{pmlim}
Let $\{\xi^n\}$ be a sequence of cadlag $\mathcal{M}_{\{h_k\}}(S\times
[0,\infty))$-valued
processes satisfying (\ref{cndrel}) for cadlag
$\mathcal{M}_{\{h_k\}}(S)$-valued processes $\{\Xi^n\}$. If the finite-dimensional
distributions of $\xi^n$ converge to the finite-dimensional
distributions of $\xi$, then the finite-dimensional
distributions of $\Xi^n$ converge to the finite-dimensional
distributions of $\Xi$ satisfying
\[
E\bigl[e^{-\int_{S\times[0,K]}f\,d\xi(t)}|\mathcal{F}_t^{\Xi}\bigr]=e^{-K\int_
S(1-e^{-f})\,d\Xi(t)}.
\]

For $n=1,2,\ldots,$ let $\xi^n$ be a cadlag $\mathcal{M}_{\{h_k\}}(S
\times[0,r_n])$-valued
process satisfying (\ref{cndrelfi}) for cadlag
$\mathcal{M}_{\{h_k\}}(S)$-valued processes $\{\Xi^n_{r_n}\}$. If $
r_n\rightarrow\infty$ and the finite-dimensional distributions of $\xi^n$ converge to the
finite-dimensional distributions of $\xi$, then the finite-dimensional
distributions of $\Xi^n_{r_n}$ converge to the finite-dimensional
distributions of $\Xi$ satisfying
\[
E\bigl[e^{-\int_{S\times[0,K]}f\,d\xi(t)}|\mathcal{F}_t^{\Xi}\bigr]=e^{-K\int_
S(1-e^{-f})\,d\Xi(t)}.
\]
\end{theorem}
\begin{pf}
Convergence of the finite-dimensional distributions
follows easily from the estimates in Lemma \ref{coxappr}.
\end{pf}

%s4.11 ###
\subsection{Martingale lemmas}

\begin{lemma}\label{mglem}
Let $\{\mathcal{F}_t\}$ and $\{\mathcal{G}_t\}$ be filtrations with
$\mathcal{G}_
t\subset\mathcal{F}_t$. Suppose
that $E[|X(t)|+\int_0^t|Y(s)|\,ds]<\infty$ for each $t$, and
\[
M(t)=X(t)-\int_0^tY(s)\,ds
\]
is an $\{\mathcal{F}_t\}$-martingale. Then
\[
\widehat{M}(t)=E[X(t)|\mathcal{G}_t]-\int_0^tE[Y(s)|\mathcal{G}_s]\,ds
\]
is a $\{\mathcal{G}_t\}$-martingale.
\end{lemma}
\begin{pf}
Let $D\in\mathcal{G}_t\subset\mathcal{F}_t$. Then
\begin{eqnarray*}
&&E\bigl[\bigl(\widehat{M}(t+r)-\widehat{M}(t)\bigr){\mathbf1}_D\bigr]\\
&&\qquad=E\biggl[\biggl(E[X(t+r)|\mathcal{G}_{t
+r}]-E[X(t)|\mathcal{G}_t]-\int_t^{t+r}E[Y(s)|\mathcal{G}_s]\,ds\biggr){\mathbf1}_
D\biggr]\\
&&\qquad=E\biggl[\biggl(X(t+r)-X(t)-\int_t^{t+r}Y(s)\,ds\biggr){\mathbf1}_D\biggr]\\
&&\qquad=0,
\end{eqnarray*}
giving the martingale property.
\end{pf}
\begin{lemma}\label{cnddc}
Let $\{\mathcal{F}_n\}$ be an increasing sequence of $\sigma$-algebras
and $
\{X_n\}$
a sequence of random variables satisfying
$E[{\sup_n}|X_n|]<\infty$ and\break $\lim_{n\rightarrow\infty}X_n=X$ a.s. Then
\[
\lim_{n\rightarrow\infty}E[X_n|\mathcal{F}_n]=E\biggl[X\Big|\bigvee_n\mathcal{F}_n
\biggr].
\]
\end{lemma}
\begin{pf}
By the martingale convergence theorem, we have
\[
E\biggl[\inf_{k\geq m}X_k\Big|\bigvee_n\mathcal{F}_n\biggr]\leq\liminf_{n\rightarrow\infty}
E[X_n|\mathcal{F}_n]\leq\limsup_{n\rightarrow\infty}E[X_n|\mathcal{F}_n
]\leq E\biggl[\sup_{k\geq m}X_k\Big|\bigvee_n\mathcal{F}_n\biggr],
\]
and the result follows by letting $m\rightarrow\infty$.
\end{pf}

%s4.12 ###
\subsection{Markov mapping theorem}\label{mpsect} The
following theorem (extending Corollary 3.5 from
\cite{K98}) plays an essential role in justifying the
particle representations and can also be used to prove
uniqueness for the corresponding measure-valued
processes. Let $(S,d)$ and $(S_0,d_0)$ be complete, separable
metric spaces, $B(S)\subset M(S)$ be the Banach space of
bounded measurable functions on $S$, with $\|f\|={\sup_{x\in S}}|f(
x)|$ and
$\overline{C}(S)\subset B(S)$ be the subspace of bounded continuous
functions. An operator $A\subset B(S)\times B(S)$ is \textit{dissipative}
if
$\Vert f_1-f_2-\varepsilon(g_1-g_2)\Vert\geq\Vert f_1-f_2\Vert$ for all $
(f_1,g_1),(f_2,g_2)\in A$ and
$\varepsilon>0$; $A$ is a \textit{pre-generator} if $A$ is dissipative and there
are sequences of functions $\mu_n\dvtx S\rightarrow\mathcal{P}(S)$ and $
\lambda_n\dvtx S\rightarrow[0,\infty)$
such that for each $(f,g)\in A$
%
%e4.26 ###
\begin{equation}\label{gencomp}
g(x)=\lim_{n\rightarrow\infty}\lambda_n(x)\int_S\bigl(
f(y)-f(x)\bigr)\mu_n(x,dy)
\end{equation}
for each $x\in S$. $A$ is \textit{graph separable}
if there exists a countable subset
$\{g_k\}\subset\mathcal{D}(A)\cap\overline{C}(S)$
such that the graph of $A$ is contained in the bounded, pointwise
closure of the linear span of $\{(g_k,Ag_k)\}$. [More precisely,
we should say that there exists $\{(g_k,h_k)\}\subset A\cap\overline{
C}(S)\times B(S)$
such that $A$ is contained in the bounded pointwise
closure of $\{(g_k,h_k)\}$, but typically $A$ is single-valued, so
we use the more intuitive notation $Ag_k$.] These two
conditions are satisfied by essentially all operators $A$
that might reasonably be thought to be generators of
Markov processes. Note that $A$ is graph separable if
$A\subset L\times L$, where $L\subset B(S)$ is separable in the sup norm
topology, for example, if $S$ is locally compact, and $L$ is
the space of continuous functions vanishing at infinity.

A collection of functions $D\subset\overline{C}(S)$ is \textit{separating} if
$\nu,\mu\in\mathcal{P}(S)$ and $\int_Sf\,d\nu=\int_Sf\,d\mu$ for all $
f\in D$ imply
$\mu=\nu$.

For an $S_0$-valued, measurable process $Y$, $\widehat{\mathcal{F}}^Y_t$
will denote the
completion of the $\sigma$-algebra\vspace*{1pt} $\sigma(Y(0),\int_0^rh(Y(s))\,d
s,r\leq t,h\in B(S_0))$.
For almost every $t$, $Y(t)$ will be $\widehat{\mathcal
{F}}_t^Y$-measurable, but in
general, $\widehat{\mathcal{F}}^Y_t$ does not contain $\mathcal
{F}^Y_t=\sigma
(Y(s)\dvtx s\leq t)$. Let
$\mathbf{T}^Y=\{t\dvtx Y(t)\mbox{ is }\widehat{\mathcal{F}}_t^Y\mbox{ measurable}
\}$. If $Y$ is cadlag and has no
fixed points of discontinuity [i.e., for every $t$,
$Y(t)=Y(t-)$ a.s.], then $\mathbf{T}^Y=[0,\infty)$. $D_S[0,\infty
)$ denotes the
space of cadlag, $S$-valued functions with the Skorohod
topology, and $M_S[0,\infty)$ denotes the space of Borel
measurable functions, $x\dvtx [0,\infty)\rightarrow S$, topologized by
convergence in Lebesgue measure.
\begin{theorem}\label{mf}Let
$(S,d)$ and $(S_0,d_0)$ be complete, separable metric spaces.
Let $A\subset\overline{C}(S)\times C(S)$ and $\psi\in C(S)$, $\psi\geq
1$. Suppose that for
each $f\in\mathcal{D}(A)$ there exists $c_f>0$ such that
%
%e4.27 ###
\begin{equation}\label{opest}
|Af(x)|\leq c_f\psi(x),\qquad  x\in A,
\end{equation}
and define $A_0f(x)=Af(x)/\psi(x)$.

Suppose that $A_0$ is a graph-separable pre-generator,
and suppose that $\mathcal{D}(A)=\mathcal{D}(A_0)$ is closed under
multiplication and
is separating. Let $\gamma\dvtx S\rightarrow S_0$ be Borel measurable, and let
$\alpha$ be a transition function from $S_0$ into $S$
[$y\in S_0\rightarrow\alpha(y,\cdot)\in\mathcal{P}(S)$ is Borel
measurable] satisfying\vspace*{1pt}
$\int h\circ\gamma(z)\alpha(y,dz)=h(y)$, $y\in S_0$, $h\in B(S_
0)$, that is,
$\alpha(y,\gamma^{-1}(y))=1$. Assume that $\widetilde{\psi}(y)\equiv
\int_S\psi(z)\alpha(y,dz)<\infty$ for each
$y\in S_0$, and define
\[
C=\biggl\{\biggl(\int_Sf(z)\alpha(\cdot,dz),\int_SAf(z)\alpha(\cdot,dz)
\biggr)\dvtx f\in\mathcal{D}(A)\biggr\}.
\]
Let $\mu_0\in\mathcal{P}(S_0)$, and define $\nu_0=\int\alpha(y,\cdot
)\mu_0(dy)$.

\begin{enumerate}[(a)]
\item[(a)] If
$\widetilde{Y}$ satisfies $\int_0^tE[\widetilde{\psi}(\widetilde{Y}(s))]\,ds
<\infty$ for all $t\geq0$, and $\widetilde{Y}$
is a solution of the martingale problem for $(C,\mu_0)$,
then there exists a solution $X$ of the martingale problem
for $(A,\nu_0)$ such that $\widetilde{Y}$ has the same distribution on
$M_{S_0}[0,\infty)$ as $Y=\gamma\circ X$. If $Y$ and $\widetilde{Y}$ are
cadlag, then $
Y$
and $\widetilde{Y}$ have the same distribution on $D_{S_0}[0,\infty
)$.

\item[(b)]For $t\in\mathbf{T}^Y$,
%
%e4.28 ###
\begin{equation}\label{rpid}
P\{X(t)\in\Gamma|\widehat{\mathcal{F}}^Y_t\}=\alpha(Y
(t),\Gamma),\qquad \Gamma\in\mathcal{B}(S).
\end{equation}

\item[(c)]If, in addition, uniqueness holds for the martingale
problem for $(A,\nu_0)$, then
uniqueness holds for the $M_{S_0}[0,\infty)$-martingale
problem for $(C,\mu_0)$. If $\widetilde{Y}$ has sample paths in $D_{
S_0}[0,\infty)$, then
uniqueness holds for the\break $D_{S_0}[0,\infty)$-martingale problem for
$(C,\mu_0)$.

\item[(d)]If uniqueness holds for the martingale problem for
$(A,\nu_0)$,
then $Y$ restricted to $\mathbf{T}^Y$ is a Markov process.
\end{enumerate}
\end{theorem}
\begin{remark}\label{discgen}
Theorem \ref{mf} can be extended to cover a large class
of generators whose range contains discontinuous
functions. (See \cite{K98}, Corollary 3.5 and Theorem
2.7.) In particular, suppose $A_1,\ldots,A_m$ satisfy the
conditions of Theorem \ref{mf} for a common domain
$\mathcal{D}=\mathcal{D}(A_1)=\cdots=\mathcal{D}(A_m)$, and $\beta
_1,\ldots
,\beta_m$ are nonnegative
functions in $M(S)$. Then the conclusions of Theorem
\ref{mf} hold for
\[
Af=\beta_1A_1f+\cdots+\beta_mA_mf.
\]
\end{remark}
\begin{pf*}{Proof of Theorem \ref{mf}}
Theorem 3.2 of \cite{K98} can be extended to operators
satisfying (\ref{opest}) by applying Corollary 1.12 of
\cite{KS01} (with the operator $B$ in that corollary set
equal zero) in place of Theorem 2.6 of \cite{K98}.
Alternatively, see Corollary 3.3 of~\cite{KN09}.
\end{pf*}

%s4.13 ###
\subsection{Uniqueness for martingale problems}\label{unmgp}
Assume that $B\subset\overline{C}(E)\times\overline{C}(E)$, that $\mathcal{D}
(B)$ is closed under
multiplication and is separating, and that existence and
uniqueness hold for the $D_E[0,\infty)$ martingale problem for $
(B,\nu)$
for each initial distribution $\nu\in\mathcal{P}(E)$.
Without loss of generality, we can assume $g\in\mathcal{D}(B)$
satisfies $0\leq g\leq1$.

By Theorem 4.10.1
of \cite{EK86}, existence and
uniqueness then follows for the $n$-particle motion
martingale problem with generator
%
%e4.29 ###
\begin{equation}\label{ngen}
B^n=\Biggl\{\Biggl(f(x),f(x)\sum_{i=1}^n\frac{Bg(x_i)}{g(x_i
)}\Biggr)\dvtx f(x)=\prod_{i=1}^ng(x_i)\Biggr\}.
\end{equation}
Actually, the cited theorem implies uniqueness for the
ordered $n$-particle motion with generator
\[
\widetilde{B}^n=\Biggl\{\Biggl(f(x),f(x)\sum_{i=1}^n\frac{Bg_i(x_i)}{g_i(x_i)}
\Biggr)\dvtx f(x)=\prod_{i=1}^ng_i(x_i),g_i\in\mathcal{D}(B)\Biggr\},
\]
but Theorem \ref{mf} can be applied to obtain uniqueness
for $B^n$ from uniqueness for~$\widetilde{B}^n$. Define $\gamma(x)
=\sum_{i=1}^n\delta_{x_i}$ and
let $\alpha(y,\cdot)$ average over all permutations of the $x_i$ in
$y=\sum_{i=1}^n\delta_{x_i}$.

Now consider a generator for a process with state space
\[
S=\bigcup_n\Biggl\{\sum_{i=1}^n\delta_{x_i}\dvtx x_i\in E\Biggr\},
\]
where we allow $n=0$, that is, no particles exist.

For $f(x,n)=\prod_{i=1}^ng(x_i)$, let
\begin{eqnarray*}
Af(x,n)&=&B^nf(x,n)\\
&&{} +f(x,n)\sum_k\lambda_k(x)\int_{E^k}\Biggl(\prod_{i=1}^
kg(z_i)-1\Biggr)\eta_k(x,dz_1,\ldots,dz_k)\\
&&{} +f(x,n)\sum_{(z_1,\ldots,z_l)\subset\{x_1,\ldots,x_n\}}
\mu(x,z_1,\ldots,z_l)\biggl(\frac1{\prod_{i=1}^lg(z_i)}-1\biggr),
\end{eqnarray*}
where $\lambda_k,\mu\geq0$ and $\eta_k$ is a transition function from $
S$ to
$E^k$. The generator has the following simple interpretation:
in between birth and death events the particles move
independently with motion determined by $B$. At rate
$\lambda_k(x)$, $k$ new particles are created with locations in $
E$
determined by $\eta_k$. At rate $\mu(x,z_1,\ldots,z_l)$, the
particles at
$z_1,\ldots,z_l$ are removed.

Let
\[
\beta(x)=\sum_k\lambda_k(x)+\sum_{(z_1,\ldots,z_l)\subset\{x_
1,\ldots,x_n\}}\mu(x,z_1,\ldots,z_l).
\]
Then for each initial distribution $\nu_0$ and each $m>0$, a localization
argument and Theorem 4.10.3 of \cite{EK86} imply
existence and uniqueness of the martingale problem for
$(A,\nu_0)$ up to the first time the solution leaves
$\{x\dvtx\beta(x)<m\}$. Consequently, existence and uniqueness hold
provided that there is a solution $X$ satisfying
$\sup_{s\leq t}\beta(X(s))<\infty$ a.s. for each $t>0$.

Essentially the same argument gives existence and
uniqueness for generators of the form (\ref{fingen})
provided $\inf_x(a(x)r-b(x))>0$ and there exists a solution
satisfying $\sup_{s\leq t}\sum a(X_i(s))<\infty$ a.s. for each $t
>0$.
\end{appendix}

% imsref loaded by lrinkeviciute, 2010-11-18 13:31:02
%

%
\printaddresses

\end{document}